\documentclass[11pt]{article}
\usepackage{graphicx,amsmath,amsfonts,amssymb,bbm}
\usepackage{psfrag}
\usepackage{epsf}
\makeatletter\@addtoreset {equation}{section}\makeatother

\newtheorem{theo}{Theorem}[section]
\newtheorem{lem}[theo]{Lemma}

\newtheorem{conj}[theo]{Conjecture}
\newtheorem{prop}[theo]{Proposition}

\newtheorem{rem}[theo]{Remark}
\newenvironment{Proof}
{\begin{trivlist} \item[]{\bf Proof. }}%
{\hspace*{\fill}$\rule{.3\baselineskip}{.35\baselineskip}$\end{trivlist}}

\newcommand{\R}{\mathbb{R}}
\newcommand{\C}{\mathbb{C}}

\newcommand{\PT}{{\cal PT}}
\newcommand{\p}{{\cal P}}
\newcommand{\T}{{\cal T}}

\renewcommand{\geq}{\geqslant}
\renewcommand{\leq}{\leqslant}

\renewcommand{\phi}{\varphi}
\newcommand{\be}{\begin{eqnarray}}
\newcommand{\ee}{\end{eqnarray}}

\newcommand{\eps}{\varepsilon}

\begin{document}

\title{\bf On the Thomas--Fermi approximation of the ground state
in a $\PT$-symmetric confining potential}

\author{Cl\'ement Gallo$^{1}$ and Dmitry Pelinovsky$^{2}$ \\
{\small $^{1}$ Institut de Math\'ematiques et de Mod\'elisation,
Universit\'e Montpellier II, 34095 Montpellier, France} \\
{\small $^{2}$ Department of Mathematics, McMaster
University, Hamilton, Ontario, Canada, L8S 4K1}  }

\date{\today}
\maketitle

\begin{abstract}
For the stationary Gross--Pitaevskii equation with harmonic real and linear imaginary potentials in the space of
one dimension, we study the ground state in the limit of large densities (large chemical potentials),
where the solution degenerates into a compact Thomas--Fermi approximation. We prove that the Thomas--Fermi
approximation can be constructed with an invertible coordinate transformation and an
unstable manifold theorem for a planar dynamical system. The
Thomas--Fermi approximation can be justified by reducing the existence problem
to the Painlev\'e-II equation, which admits a unique global Hastings--McLeod solution.
We illustrate numerically that an iterative approach to solving the existence problem
converges but give no analytical proof of this result.
Generalizations are discussed for the stationary Gross--Pitaevskii equation
with harmonic real and localized imaginary potentials.
\end{abstract}

\section{Introduction}

Ground states of the repulsive Bose--Einstein condensates placed in a harmonic (magnetic)
confinement are global minimizers of the Gross--Pitaevskii energy \cite{A}. For the large--density atomic gas,
these ground states are well approximated by a compact function, which is referred to as
the Thomas--Fermi approximation. The Thomas--Fermi approximation was rigorously
justified using calculus of variations \cite{IM}.

It was also discovered in several independent studies \cite{ADP,BTNN,KK} that the
nearly compact  Thomas--Fermi approximation of the ground state has a
super-exponential spatial decay outside a transitional layer, where the ground state
satisfies the Painlev\'e-II equation \cite{FIKN}. A particular solution
of the Painlev\'e-II equation referred to as the Hastings-McLeod solution \cite{HM,HM2} reconstructs
the Thomas--Fermi approximation by means of a change of dependent and independent variables \cite[Section 1.2.3]{A}.

Justification of the Painlev\'e-II equation in the context of the radially symmetrical ground states
was developed in our previous work \cite{GP}, where it was shown that the Painlev\'e-II equation remains valid
uniformly on the spatial scale if the confining potential is purely harmonic. These results were further used for
several purposes. Expansions of energy for the Thomas--Fermi approximation were studied in \cite{G}. Excited states
of the stationary Gross--Pitaevskii equation were constructed with the method of Lyapunov--Schmidt reductions
in \cite{P}. More general non-radial trapping potentials were included by Karali and Sourdis \cite{KS},
where the Painlev\'e-II equation does not hold uniformly on the spatial scale but is valid nevertheless
in a transitional layer.

Very recently, ground states of the repulsive Bose--Einstein condensates were considered
under the presence of a complex-valued potential, which expresses a certain balance between
losses and gains occurring in the atomic gases. This potential is symmetric with respect to the simultaneous
parity ($\p$) and time-reversal ($\T$) transformations, hence it is referred to as the
$\PT$-symmetric potential. Thomas--Fermi approximations in a
localized $\PT$-symmetric potential added to the harmonic potential were numerically considered
in \cite{AKFC}. Ground and excited states in a linear $\PT$-symmetric potential added to the harmonic
potential were numerically constructed in \cite{ZK}. In both works, it was discovered that
the existence of the ground state may be fragile in the presence of the $\PT$-symmetric potential,
because the branch of the ground state may disappear due to the coalescence with the branch of the first
excited state. Precise predictions on where this saddle-node bifurcation occurs and whether
the ground state can be extended to the Thomas--Fermi (large--density) limit were not detailed in these works.

In the present work, we study existence of the ground state in the Thomas--Fermi limit
for the Gross--Pitaevskii equation with a $\PT$-symmetric potential. Because no variational
principle can be formulated for the $\PT$-symmetric potential, existence of ground states
can not be established using calculus of variations. We use again the transformation
of the stationary Gross--Pitaevskii equation to the Painlev\'e-II equation with the Hastings-McLeod solution.
Persistence of the Hastings-McLeod solution is analyzed analytically and numerically
with an iterative approach.

Our starting point is the stationary Gross--Pitaevskii equation with the $\PT$-symmetric harmonic potential
\begin{equation}
\label{stationaryGP1}
\mu U(X) = \left( -\partial_X^2 + X^2 + 2 i \alpha X + |U(X)|^2\right) U(X),  \quad
X \in \R,
\end{equation}
where $\mu \in \mathbb{R}$ is the chemical potential, $\alpha \in \mathbb{R}$ is the gain--loss coefficient,
and $U : \R \to \C$ is the wave function for the steady state. The spectrum of the linearized operator
$$
L_0 := -\partial_X^2 + X^2 + 2 i \alpha X = -\partial_X^2 + (X + i \alpha)^2 + \alpha^2
$$
coincides with the spectrum of the operator $\tilde{L}_0 = -\partial_Z^2 + Z^2 + \alpha^2$. Therefore,
the spectrum of $L_0$ is purely discrete, real, and bounded from below.
The ground state of the stationary Gross--Pitaevskii equation (\ref{stationaryGP1})
bifurcates from the smallest eigenvalue
$\mu_0 = 1 + \alpha^2$ of the operator $L_0$ and exists for $\mu \geq \mu_0$. This local bifurcation of the ground state (as well as
excited states) was formally considered by Zezyulin and Konotop \cite{ZK}.

As the Thomas--Fermi approximation is derived in the large-density limit $\mu \to \infty$,
we introduce the change of variables $\mu = \eps^{-1}$ and $U(X) = \eps^{-1/2} u(\eps^{1/2} X)$
with small positive $\eps$. The stationary Gross--Pitaevskii equation (\ref{stationaryGP1}) is now
written in the form
\begin{equation}
\label{stationaryGP2}
\left( 1 - x^2 - 2 i \alpha \eps^{1/2} x - |u(x)|^2 \right) u(x) = - \eps^2 u''(x),  \quad
x \in \R,
\end{equation}
where $x = \eps^{1/2} X$. To incorporate the complex phase of $u$ produced by
the gain-loss term, we use the polar form $u = \varphi e^{i \theta}$ for the ground state
solutions with $|u(x)| > 0$ for all $x \in \mathbb{R}$. Splitting
the scalar equation (\ref{stationaryGP2}) for real and imaginary parts, we obtain the system
\begin{equation}
\label{stationaryGP3}
\left\{ \begin{array}{l}
\left( 1 - x^2 - \varphi^2(x) - \eps^2 (\theta')^2(x) \right) \varphi(x) = - \eps^2 \varphi''(x),  \\
\left( \varphi^2 \theta' \right)'(x) = 2 \alpha \eps^{-3/2} x \varphi^2(x), \end{array} \right. \quad x \in \R.
\end{equation}
As we are looking for spatially decaying solutions with $\varphi^2 \theta'(x) \to 0$ as $|x| \to \infty$,
it is clear that the following constraint must be satisfied:
\begin{equation}
\int_{\mathbb{R}} x \varphi^2(x) dx = 0.
\end{equation}
In particular, if $\varphi^2$ and $\theta'$ are even in $x$, the constraint is satisfied, and
the parity requirement implies that the stationary solution is $\PT$-symmetric with
$u(-x) = \bar{u}(x) e^{i \theta_0}$, where $\theta_0 \in \mathbb{R}$.
Note in passing that the question whether $\PT$-symmetric equations may
admit non-$\PT$-symmetric spatially decaying solutions is open, recent results in
this direction were obtained by Yang \cite{Yang} with perturbation techniques.

Using the scaled variable $\xi(x) := \eps \theta'(x)$ and scaled
parameter $\alpha = \eps^{1/2} \eta$, we obtain the final form
of the existence problem:
\begin{equation}
\label{stationaryGP}
\left\{ \begin{array}{l}
\left( 1 - x^2 - \varphi^2(x) - \xi^2(x) \right) \varphi(x) = - \eps^2 \varphi''(x),  \\
\left( \varphi^2 \xi \right)'(x) = 2 \eta x \varphi^2(x), \end{array} \right. \quad x \in \R.
\end{equation}
The existence problem has two parameters $\eta$ and $\eps$ and we are looking
for the ground state with even and strictly positive $\varphi$ in the limit of small $\eps$.
The other parameter $\eta$ can be either $\eps$-independent or $\eps$-dependent
and we shall later specify the conditions
on this parameter to ensure that the ground state exists in the limit of small $\eps$.

When $\eta = 0$, we can choose $\xi \equiv 0$ and the existence problem (\ref{stationaryGP})
reduces to the stationary Gross--Pitaevskii equation, studied in our previous work \cite{GP}.
In the most general case, we can solve the second equation of system (\ref{stationaryGP}) uniquely
from the boundary condition $\lim\limits_{x \to \pm \infty} \varphi^2(x) \xi(x) = 0$. In this way, we
obtain the integral representation
\begin{equation}
\label{component-xi}
\xi(x) = \frac{2 \eta}{\varphi^2(x)} \int_{-\infty}^x s \varphi^2(s) ds,
\end{equation}
which allows us to close the first equation of system
(\ref{stationaryGP}) as an integro--differential equation.

The formal Thomas--Fermi limit corresponds to the solution of the truncated problem
\begin{equation}
\label{TFlimit}
\left\{ \begin{array}{ll} 1 - x^2 - \varphi^2(x) - \xi^2(x) = 0,  \\
\left( \varphi^2 \xi \right)'(x) = 2 \eta x \varphi^2(x), \end{array}\right. \quad x \in [-1,1],
\end{equation}
subject to the boundary conditions $\varphi(\pm 1) = \xi(\pm 1) = 0$.
In the following theorem, we state the existence of suitable solutions to the limiting
problem (\ref{TFlimit}) for a sufficiently small but $\eps$-independent $\eta$.
Since the component $\xi$ is uniquely determined by (\ref{component-xi}),
we set
\begin{equation}
\label{component-xi-truncated}
\xi(x) = \frac{2 \eta}{\varphi^2(x)} \int_{-1}^x s \varphi^2(s) ds, \quad x \in (-1,1),
\end{equation}
and state the result in terms of $\varphi$ only.

\begin{theo}
There exists $\eta_0>0$ such that for any $|\eta|< \eta_0$,
the truncated existence problem (\ref{TFlimit})--(\ref{component-xi-truncated})
admits a unique solution $\varphi_{\rm TF} \in C^{\infty}(-1,1)$ such that
$\varphi_{\rm TF}(x) > 0$ for all $x \in (-1,1)$ and
\begin{equation}
\label{boundary-1}
\varphi^2_{\rm TF}(x) = 1 - x^2 + \mathcal{O}((1-x^2)^2) \quad \mbox{\rm as} \;\; |x| \to 1.
\end{equation}
\label{theorem-1}
\end{theo}

\begin{figure}[h]
\begin{center}
\includegraphics[scale=0.35]{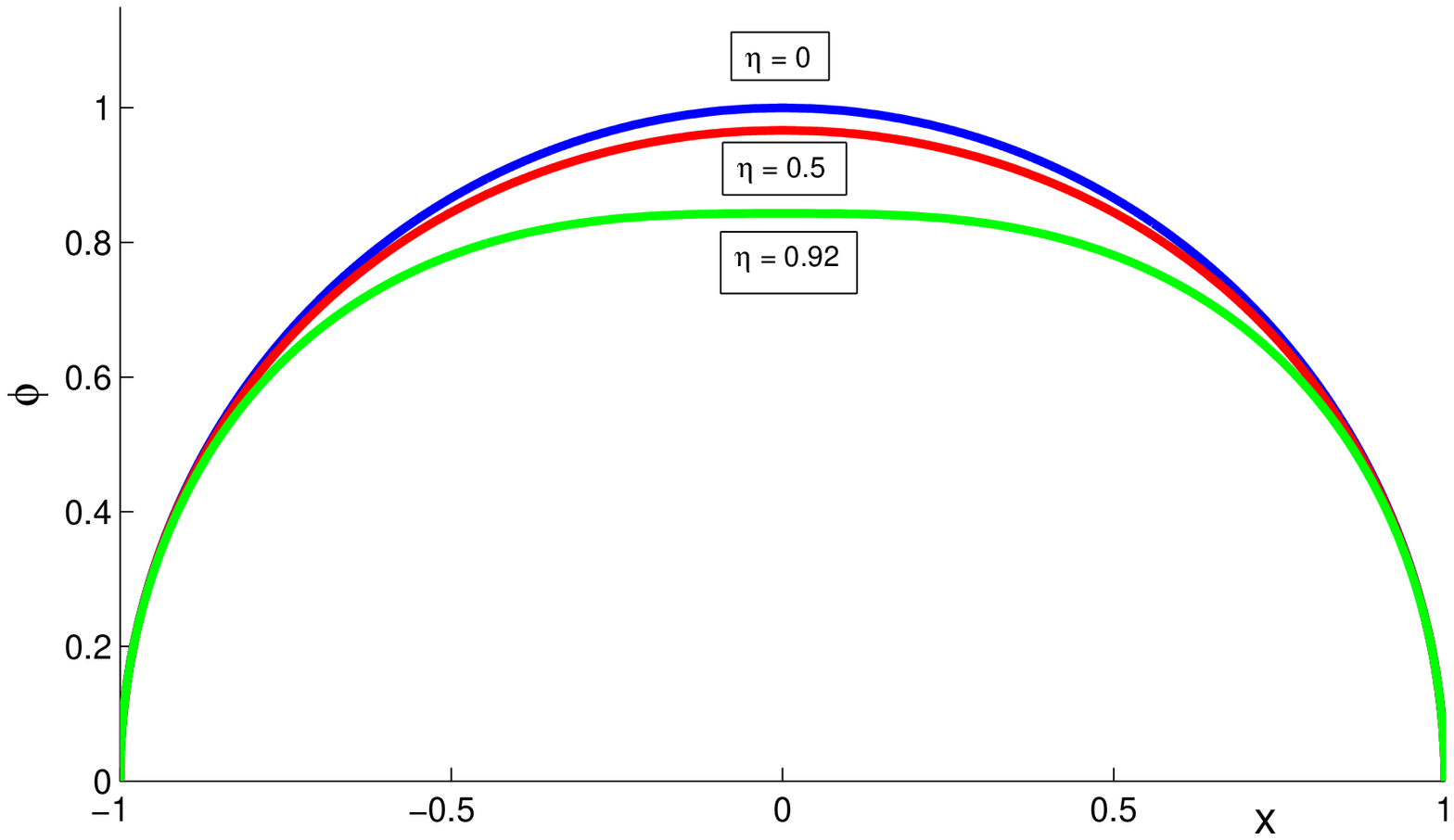}
\includegraphics[scale=0.3]{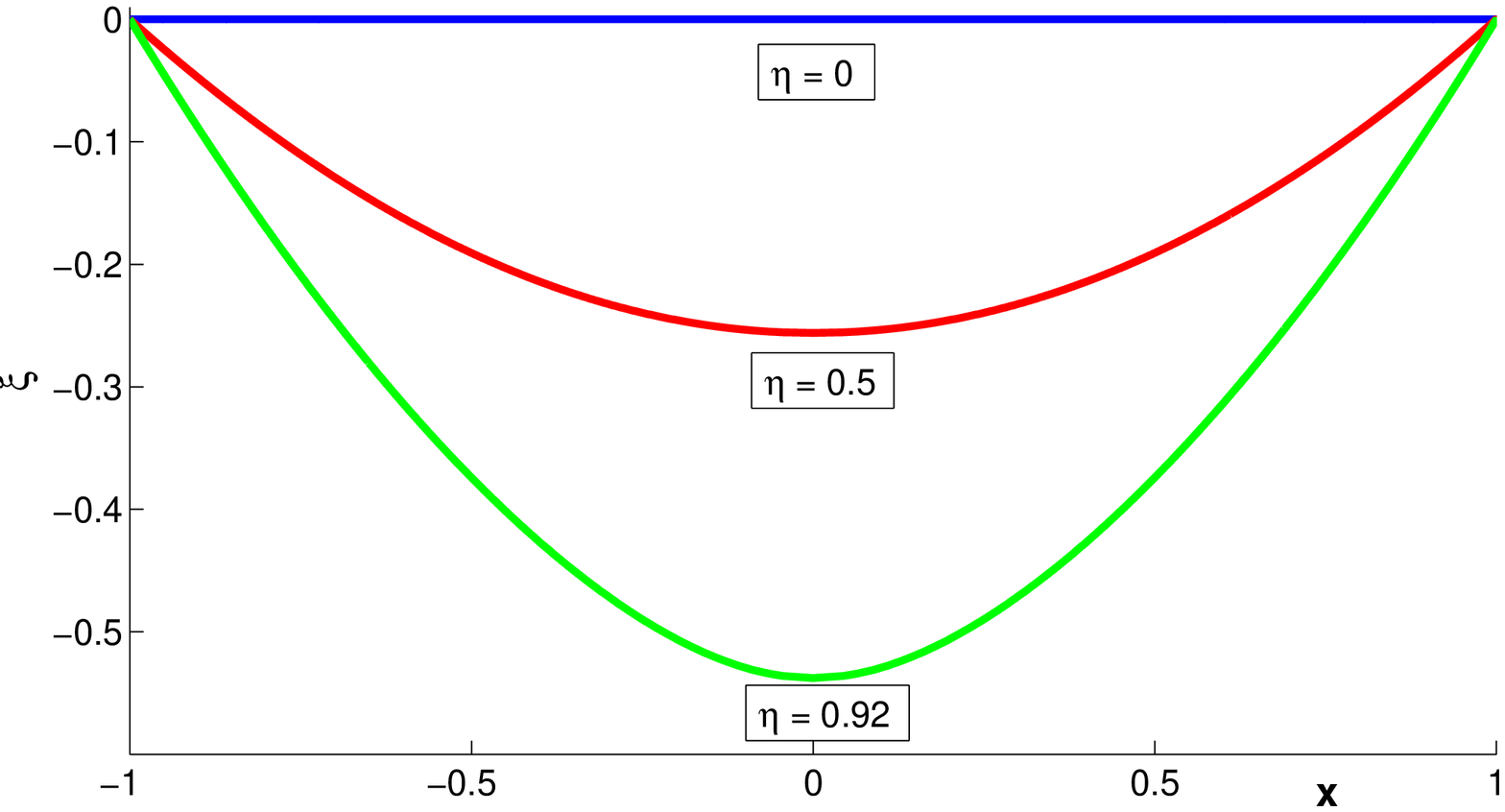}
\end{center}
\caption{\label{fig-TFprofile}
Components $\varphi$ (left) and $\xi$ (right) for the numerical solution
to the limiting problem (\ref{TFlimit}) for three different values of $\eta$.}
\end{figure}

Figure \ref{fig-TFprofile} illustrates components $\varphi$ (left) and $\xi$ (right)
of the Thomas--Fermi solution in Theorem \ref{theorem-1} for three different values of $\eta$.
The numerical solution is obtained with the fourth-order Runge--Kutta method applied
to the closed first-order differential equation for variable $\xi$,
after the variable $\varphi^2$ is eliminated from the system (\ref{TFlimit}).
The solution terminates at $\eta_0 \approx 0.93$ because the derivative of $\xi$ diverges near $x = 0$.
The numerical approximations illustrate the statement of Theorem \ref{theorem-1} that
the Thomas--Fermi approximation exists only for $|\eta| < \eta_0$, where the value of $\eta_0$ is finite.

\begin{rem}
From Theorem \ref{theorem-1} and numerical illustrations, we can see that
the Thomas--Fermi radius ($|x| = 1$ in this particular case) is independent of the gain-loss parameter $\eta$
and that the $\PT$-symmetric linear potential leads to the decrease of the ground state amplitude
$\varphi$ near the center of the harmonic potential ($x = 0$). These two facts
appear to be universal for spatially decaying $\PT$-symmetric potentials.
\label{remark-1}
\end{rem}

The limit leading to the Painlev\'e-II equation appears after the formal
change of dependent and independent variables near the Thomas--Fermi radius ($x = 1$):
\begin{equation}
\label{asymptotic-transformation}
\varphi(x) = \eps^{1/3} \nu(y), \quad \xi(x) = \eps^{2/3} \chi(y), \quad y = \frac{1-x^2}{\eps^{2/3}}.
\end{equation}
The new variables satisfy the modified existence problem:
\begin{equation}
\label{P2limit}
\left\{ \begin{array}{l}
4 \nu''(y) + y \nu(y) - \nu^3(y) = \eps^{2/3} \left( 4 y \nu''(y) + 2 \nu'(y) + \chi^2(y) \nu(y) \right),  \\
\left( \nu^2 \chi \right)'(y) = -\eta \nu^2(y), \end{array} \right. \quad y \in (-\infty,\eps^{-2/3}),
\end{equation}
subject to the decay condition $\nu(y) \to 0$ as $y \to -\infty$.
The truncated version of the first equation in system (\ref{P2limit})
is the Painlev\'e-II equation
\begin{equation}
\label{P2-eq} 4 \nu''(y) + y \nu(y) - \nu^3(y) = 0,\quad
y\in\R,
\end{equation}
which admits a unique solution $\nu_0$ \cite{HM}
satisfying the following asymptotic behavior \cite{FIKN}
\begin{equation}
\label{P-bc}
\nu_0(y) = \left\{ \begin{array}{l}
y^{1/2} - \frac{1}{2}  y^{-5/2} + \mathcal{O}(y^{-11/2}) \quad \quad \quad \quad \quad \quad \mbox{as} \quad \; y\to +\infty, \\
\pi^{-1/2} |y|^{-1/4} e^{-\frac{1}{3} |y|^{3/2}} \left(1 + \mathcal{O}(|y|^{-3/4}) \right) \quad \; \mbox{as} \quad y\to-\infty. \end{array} \right.
\end{equation}
This function $\nu_0$ is referred to as the Hastings--McLeod
solution of the Painlev\'e-II equation (\ref{P2-eq}). Moreover, the asymptotic expansion as $y\to +\infty$ in (\ref{P-bc})
can be differentiated term by term.

The persistence of the Hastings--McLeod solution $\nu_0$ with respect to
small perturbation terms in $\eps$ needs to be considered within the modified existence problem (\ref{P2limit}).
For technical reasons, it is easier to work with a small $\eps$-dependent $\eta$ but
even with this simplification, we obtain a partial progress towards the proof of
persistence. Since the component $\chi$ is uniquely determined
by integrating the second equation of system (\ref{P2limit}), we set
\begin{equation}
\label{chi}
\chi(y) = -\frac{\eta}{\nu^2(y)} \int_{-\infty}^y \nu^2(s) ds
\end{equation}
and state the desired result in terms of $\nu$ only.

\begin{conj}
\label{theorem-2}
Let $\nu_0$ be the Hastings--McLeod solution of the Painlev\'e-II
equation, defined by (\ref{P2-eq}) and (\ref{P-bc}). For any $q > \frac{5}{6}$,
there exist $\eps_q > 0$, $\eta_q > 0$, and $C_q > 0$ such that for every $\eps \in (0,\eps_q)$
and $|\eta| < \eta_q \eps^{q}$, the coupled system (\ref{P2limit}) admits a unique solution
$\nu_{\rm P} \in C^{\infty}(-\infty,\eps^{-2/3})$ such that
$\nu_{\rm P}(y) > 0$ for all $y \in (-\infty,\eps^{-2/3})$,
$\nu_{\rm P}(y) \to 0$ as $y \to -\infty$, and
\begin{equation}
\label{upper-bound-theorem}
\sup_{y \in (-\infty,\eps^{-2/3})} \left| \nu_{\rm P}(y) - \nu_0(y) \right| \leq
C_q \left\{ \begin{array}{l} \eps^{2q - 4/3} \left| \log(\eps) \right|^{1/2}, \quad q \leq 1, \\
\eps^{2/3}, \quad  \quad  \quad  \quad  \quad  \quad  \quad q > 1. \end{array} \right.
\end{equation}
\end{conj}

\begin{rem}
It follows from the bound (\ref{upper-bound-theorem}) that for every $x \in (-1,1)$, we have
\begin{equation}
\label{boundary-2}
\eps^{2/3} \nu^2_{\rm P}(y) \to 1 - x^2 \quad \mbox{\rm as} \quad \eps \to 0,
\end{equation}
where $y = \frac{1-x^2}{\eps^{2/3}}$.
\end{rem}

The bound (\ref{upper-bound-theorem}) in Conjecture \ref{theorem-2} is found from the rigorous analysis
of the solution of the persistence problem if a suitable bounded
function $\chi$ is substituted in the first equation of the system (\ref{P2limit}).
For this reduced problem, we can prove existence of the solution for the component $\nu$
satisfying the bound (\ref{upper-bound-theorem}) (see Theorem \ref{theorem-map-chi-to-nu} below).
When this solution for
the component $\nu$ is used in the integral equation (\ref{chi}), we can
also fully characterize properties of the component $\chi$ (see Lemmas \ref{theorem-map-nu-to-chi}
and \ref{lemma-mapping} below). By alternating
solutions of these two uncoupled problems, we can develop a simple iterative method,
which approximates numerically solutions of the coupled system (\ref{P2limit}). Although
this numerical method is found to converge extremely fast, we still lack nice Lipschitz
properties of the integral equation (\ref{chi}) in order to achieve a rigorous proof
of the statement in Conjecture \ref{theorem-2}.

\begin{figure}[h]
\begin{center}
\includegraphics[scale=0.33]{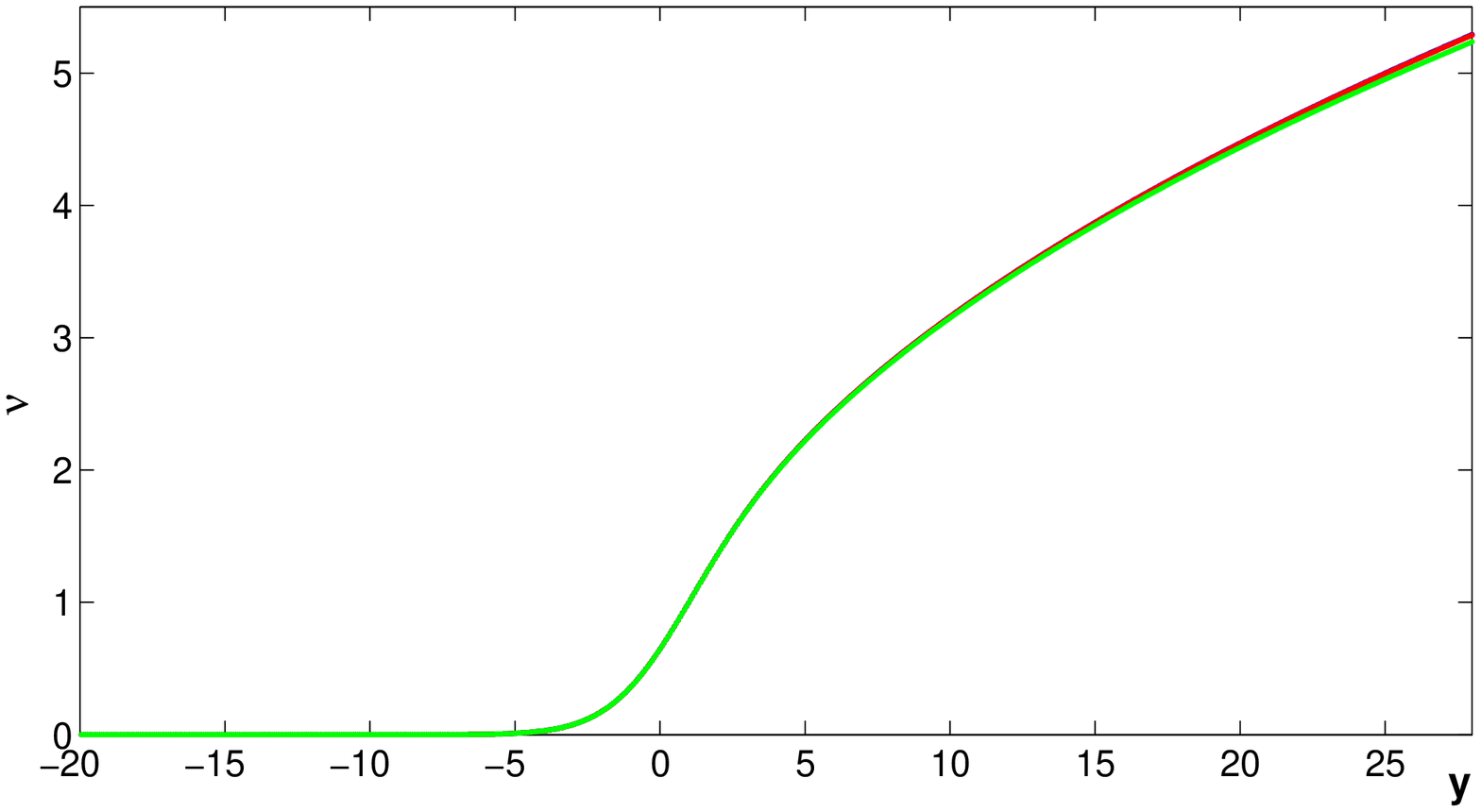}
\includegraphics[scale=0.32]{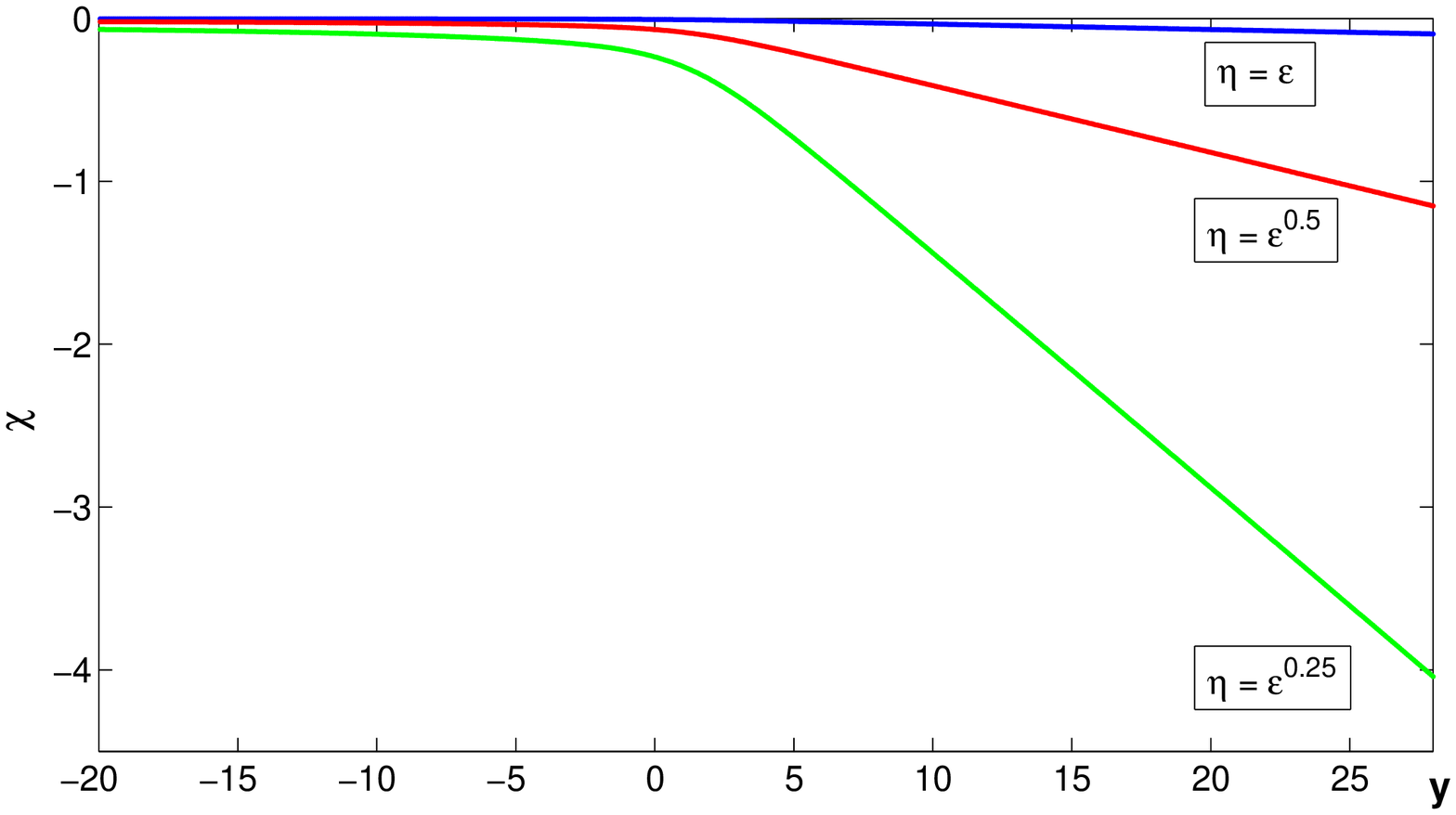}
\end{center}
\caption{\label{fig-HMprofile}
Components $\nu$ (left) and $\chi$ (right) for the numerical solution
to the coupled system (\ref{P2limit}) with $\eps = 0.0067$ and three
different values of $\eta$.}
\end{figure}

Figure \ref{fig-HMprofile} illustrates components $\nu$ (left) and $\chi$ (right)
of the numerical solution of the coupled system (\ref{P2limit}) for $\eps = 0.0067$
and three different values of $\eta$. The numerical solution is obtained with
the iterative method described above. Because the component $\nu$ is close to the
Hastings--McLeod solution $\nu_0$, the difference between the three cases of $\eta$
is not visible on the left panel of the figure. The convergence of the numerical
method is lost for $\eta \approx \eps^{0.15}$,  which may signal that no solution
of the coupled system (\ref{P2limit}) exists for such large values of $\eta$.

This paper is organized as follows. Section 2 gives the proof of Theorem \ref{theorem-1}
on the existence of the compact Thomas--Fermi approximation. Section 3 describes properties of
decaying solutions of the system (\ref{stationaryGP}) closed with the integral equation (\ref{component-xi}).
Section 4 gives the details of how a solution of the persistence problem can be obtained if a suitable bounded
function $\chi$ is substituted in the first equation of the system (\ref{P2limit})
without computing it from the integral equation (\ref{chi}). Section 5 is devoted
to the study of the integral equation (\ref{chi}). Section 6 illustrates numerically the
persistence of the Hastings--McLeod solution beyond the Painlev\'e-II equation in Conjecture \ref{theorem-2}.
Section 7 discusses generalizations of our results to spatially decaying $\PT$-symmetric potentials
superposed with the harmonic confining potential.

\section{Proof of Theorem \ref{theorem-1}}

First, we rewrite the truncated problem (\ref{TFlimit}) in terms of
the variable $z := 1-x^2$ and introduce a new function $\omega(z) := \varphi^2(x)$.
Hence the truncated problem (\ref{TFlimit}) is rewritten as
\begin{equation}
\label{TFlimit-new}
\left\{ \begin{array}{ll} \omega(z) = z - \xi^2(z), \\
\left( \omega \xi \right)'(z) = - \eta \omega(z), \end{array}\right.  \quad z \in [0,1],
\end{equation}
subject to the boundary conditions $\omega(0) = \xi(0) = 0$. Here and in what follows,
we use the same notation $\xi$ for the function of variables $x$ and
$z$. If $(\omega,\xi)\in \mathcal{C}^1([0,1]) \times \mathcal{C}^1([0,1]) $ solves
(\ref{TFlimit-new}) and
$$
\omega'(z)=1 - 2 \xi(z) \xi'(z) > 0, \quad z \in [0,1],
$$
then $\omega$ is an isomorphism between $[0,1]$ and $[0,\omega(1)]$,
where $\omega(1)=1 - \xi^2(1) \in (0,1)$.

In terms of the new variable $Z=\omega(z)$, the system of two
equations (\ref{TFlimit-new}) can be closed as a first-order
non-autonomous differential equation. By
the chain rule, this equation is
\begin{equation}
\label{eq-xi-Z}
\frac{d}{dZ} (Z \xi(Z)) = - \eta Z \left( 1 + 2
\xi(Z) \xi'(Z)\right), \quad Z \in [0,\omega(1)],
\end{equation}
starting with $\xi(0)=0$, where we have again used the same notation
$\xi$ for the function of variables $z$ and $Z$.
Formally, equation (\ref{eq-xi-Z}) is solved by the power series expansion given by
\begin{equation}
\label{power-series-expansion}
\xi(Z) = -\frac{1}{2} \eta Z \left[ 1 + \frac{1}{3} \eta^2 Z + \frac{1}{4} \eta^4 Z^2 + \mathcal{O}(\eta^6 Z^3) \right].
\end{equation}
This expansion suggests us to write $\xi$ under the form
\begin{equation}
\label{variables-xi-psi}
\xi(Z) = -\frac{1}{2} \eta Z \psi(\zeta), \quad \zeta := \eta^2 Z.
\end{equation}
Note that since $\xi\in\mathcal{C}^1([0,\omega(1)])$ satisfies the power series
expansion (\ref{power-series-expansion}), $\psi$ satisfies the boundary condition $\psi(0) =
1$. Moreover, straightforward substitution imply that $\psi$ solves the first-order
differential equation
\begin{equation}
\label{first-order-non}
\frac{d \psi}{d \zeta} = \frac{2(1 - \psi) + \zeta \psi^2}{\zeta (1 -
  \zeta \psi)},\quad \zeta\in [0,\eta^2\omega(1)],
\end{equation}
starting with $\psi(0) = 1$.

In order to prove Theorem \ref{theorem-1}, we first prove existence
and uniqueness of a solution for $\psi(\zeta)$.
For this purpose, we shall transform the first-order non-autonomous equation (\ref{first-order-non}) into a planar
dynamical system, where the point $(\zeta,\psi) = (0,1)$ is an equilibrium state with a unique unstable
manifold extending to the domain $\zeta > 0$. To do this, we set $\tau := \log(\zeta)$ as an
evolutionary variable of the planar dynamical system and rewrite equation (\ref{first-order-non})
in the dynamical system form
\begin{equation}
\label{second-order}
\left\{ \begin{array}{l} \dot{\zeta} = \zeta, \\
\dot{\psi} = \frac{2(1 - \psi) + \zeta \psi^2}{1 - \zeta \psi}, \end{array} \right.
\end{equation}
where the dot stands for the derivative in $\tau$. We can see that $(\zeta,\psi) = (0,1)$
is a saddle point of the dynamical system (\ref{second-order}) and that the dynamical system
is analytic near this point.

The stable manifold of the linearized system at the critical point $(0,1)$ corresponds to the eigenvalue $-2$
and is the line $\zeta = 0$. The unstable manifold
of the linearized system at the critical point $(0,1)$
corresponds to the eigenvalue $\lambda = 1$ and is the line $\psi - 1 = \frac{1}{3} \zeta$,
which also follows from the power series expansion (\ref{power-series-expansion}).
There exists a unique solution of the linearized system for $\zeta > 0$
such that $\psi(\tau) \to 1$ and $\zeta(\tau) \to 0$ as $\tau \to -\infty$.
By the Unstable Manifold Theorem, there exists a unique solution of
the full nonlinear system (\ref{second-order}) with the same properties
and this solution is tangent to the unstable manifold of the linearized system at $(0,1)$
in the sense that
$$
\lim_{\zeta \to 0} \frac{\psi - 1 - \frac{1}{3} \zeta}{\zeta} = 0,
$$
again in agreement with the power series expansion (\ref{power-series-expansion}).
This solution exists at least locally, e.g. for $(-\infty,\tau_0)$ for
some $\tau_0 \in \mathbb{R}$. It is not clear if it exists globally or not,
because the unstable manifold on the plane $(\zeta,\psi)$ may intersect
the curve $\zeta \psi = 1$, where the dynamical system (\ref{second-order}) is singular.

We transfer now the result of the Unstable Manifold Theorem back to the solutions
of the truncated problem (\ref{TFlimit-new}). The solution to the
dynamical system (\ref{second-order}) we have constructed for $\tau\in
(-\infty,\tau_0)$ provides a solution $\psi(\zeta)$ to the first-order equation
(\ref{first-order-non}) for $\zeta \in [0,\zeta_0)$, with $\zeta_0 := e^{\tau_0}$.

Using the scaling transformation (\ref{variables-xi-psi}), we obtain
the existence of a solution $\xi(Z)$ to the first-order equation
(\ref{eq-xi-Z}) for $Z \in [0,\zeta_0 \eta^{-2})$.
For $|\eta|<\eta_0=\zeta_0^{1/2}$, this interval includes the interval $[0,1]$.
It follows from (\ref{variables-xi-psi}) and (\ref{first-order-non})
that there is a positive constant $C$ such that
\begin{equation}
|\xi(Z)| + |\xi'(Z)| \leq C |\eta| \quad \mbox{\rm for all} \; Z \in [0,1].
\end{equation}
Therefore, for sufficiently small values of $\eta$, the map $\alpha:
Z\rightarrow Z + \xi^2(Z)$ is an isomorphism from $[0,1]$ to
$[0,z_0]$, where $z_0=1+\xi(1)^2>1$. Indeed, if $\eta$ is small, then
$$
\alpha'(Z)=1 + 2 \xi(Z) \xi'(Z) > 0  \quad \mbox{\rm for all} \; Z \in [0,1].
$$
Then, $\omega(z)=z-\xi(\alpha^{-1}(z))^2$ and
$\xi(z)=\xi(\alpha^{-1}(z))$ define a smooth solution of the truncated
problem (\ref{TFlimit-new}) for all $z \in [0,1]$.
 By construction, the solution satisfies
$\omega(z) > 0$ for all $z \in (0,1]$ and $\omega(z) = z +
\mathcal{O}(z^2)$ as $z \to 0$. Defining
$\phi_{\rm TF}(x)=\sqrt{\omega(1-x^2)}$ for $x\in [-1,1]$,
we complete the proof of Theorem \ref{theorem-1}.

\section{Properties of decaying solutions}

Here we assume the existence of an even, $\mathcal{C}^1$, spatially decaying solution of the system of
differential equations (\ref{stationaryGP}). We prove that such a solution has a fast decay at infinity
with a specific rate and remains positive at least outside the Thomas--Fermi interval $[-1,1]$.
Both parameters $\eps$ and $\eta$ are considered to be positive and fixed.

\begin{lem}
Assume that $\varphi$ is an even $\mathcal{C}^1$ solution of system (\ref{stationaryGP})
such that $\varphi(x) \to 0$ as $|x| \to \infty$ and satisfies for large values of $|x|$,
\begin{equation}
\label{constraint-1-prep}
\varphi(x) > 0, \quad \left| \int_{-\infty}^x s \varphi^2(s) ds \right| < \infty.
\end{equation}
Assume that $1-\eps - \eps^2 \eta^2 > 0$. Then, there is $\gamma > 0$ such that
\begin{equation}
\label{constraint-1}
\varphi(x) \underset{|x| \to \infty}{\sim} \gamma |x|^{\frac{1-\eps - \eps^2 \eta^2}{2\eps}} e^{-\frac{x^2}{2\eps}},
\end{equation}
and $\varphi(x) > 0$ for all $|x| \geq 1$.
\label{lemma-ode}
\end{lem}

\begin{Proof}
We justify the decay (\ref{constraint-1}) from the Unstable/Stable Manifold Theorem
and the WKB theory. Let us consider decaying solutions of the linear second-order differential equation
\begin{equation}
\label{linear-GP}
\eps^2 \varphi_{\infty}''(x) + \left( 1 - x^2 - \xi(x)^2 \right) \varphi_{\infty}(x) = 0,
\end{equation}
where $\xi$ is defined by the integral formula
(\ref{component-xi}) computed at $\varphi = \varphi_{\infty}$.
By the WKB method without turning points \cite[Chapter 7.2]{Miller},
for a fixed positive $\eps$, decaying solutions of (\ref{linear-GP}) satisfying (\ref{constraint-1-prep})
exist and are all proportional to the particular solution $\phi_\infty$ given by
\begin{equation}
\label{expression-A-B}
\varphi_{\infty}(x) = \frac{e^{\eps^{-1} \int B(x) dx}}{\sqrt{B(x)}} ,
\end{equation}
where $B(x)$ satisfies
\begin{eqnarray}
\label{expression-B1}
B(x) = \sqrt{x^2 + \xi(x)^2 -1 + \frac{\eps^2 (2 B B'' - 3 (B')^2)}{4 B^2}}.
\end{eqnarray}
Since $B(x) \to +\infty$ as $x \to -\infty$, integration by parts yields
\begin{eqnarray*}
\int_{-\infty}^x s \varphi^2_{\infty}(s) ds = - \frac{\eps |x|}{2 B(x)^2}
e^{2\eps^{-1} \int B(x) dx} \left[ 1 + \mathcal{O}\left(\frac{\eps}{|x| B(x)}\right)
+ \mathcal{O}\left(\frac{\eps B'(x)}{B^2(x)}\right)  \right]\quad {\rm as\ }x\to -\infty.
\end{eqnarray*}
It follows from the integral formula (\ref{component-xi}) that
\begin{eqnarray}
\label{expression-B2}
\xi(x) = -\frac{\eps \eta |x|}{B(x)} \left[ 1 + \mathcal{O}\left(\frac{\eps}{|x| B(x)}\right)
+ \mathcal{O}\left(\frac{\eps B'(x)}{B^2(x)}\right)  \right]\quad {\rm as\ }x\to -\infty.
\end{eqnarray}
From (\ref{expression-B1}) and (\ref{expression-B2}),
we deduce that $\xi(x)=o(|x|)$ and $B(x) = \mathcal{O}(x)$
as $x \to -\infty$. By Taylor expansions, this further specifies the
asymptotic expansions
$$
\xi(x) = - \eps \eta \left[ 1 + \mathcal{O}\left(\frac{1}{x^2}\right)\right] \quad {\rm as\ }x\to -\infty
$$
and
$$
B(x)= |x| + \frac{\eps^2 \eta^2 - 1}{2 |x|} + \mathcal{O}\left(\frac{1}{|x|^3}\right) \quad {\rm as\ }x\to -\infty.
$$
By using (\ref{expression-A-B}), we obtain
\begin{equation*}
\varphi_{\infty}(x) = |x|^{\frac{1-\eps-\eps^2 \eta^2}{2 \eps}} e^{-\frac{x^2}{2\eps}}
\left[ 1 + \mathcal{O}\left(\frac{1}{x^2}\right) \right] \quad \mbox{\rm as} \quad x \to -\infty.
\end{equation*}
This asymptotic decay recovers (\ref{constraint-1}) by the Unstable Manifold
Theorem, which states that decaying solutions $\varphi$ of system (\ref{stationaryGP})
are all proportional to the decaying solution $\varphi_{\infty}$
of the linear equation (\ref{linear-GP}) as $x\to -\infty$.

To justify positivity of $\phi$, we represent the first equation in (\ref{stationaryGP}) as follows:
\begin{equation}
\label{eq-second-order}
\eps^2 \frac{d^2 \varphi}{d x^2} = (x^2 - 1 + \xi^2 + \varphi^2) \varphi.
\end{equation}
By the decay (\ref{constraint-1}), we have $\varphi(x) > 0$ and $\varphi'(x) > 0$ for large negative values of $x$.
Then $\varphi''(x) > 0$ for all $x \in (-\infty,-1)$, so that $\varphi(x) > 0$
for all $x \in (-\infty,-1]$.
\end{Proof}

\section{Mapping $\chi \to \nu$}

Here we consider system (\ref{P2limit}) for a family of functions $\chi \in L^{\infty}(-\infty,\eps^{-2/3})$, which
depend on $\eps$ and $\eta$. We assume that there are constants $C_+ \geq 1$ and $C_- > 0$ such that
for every $\eps > 0$ small enough and every $\eta \in \mathbb{R}$, the function $\chi$ satisfies
\begin{eqnarray}
\label{assumption-chi-1}
\frac{1}{C_+} |\eta| y \leq |\chi(y)| \leq C_+ |\eta| (1+y), \quad y \in (0,\eps^{-2/3})
\end{eqnarray}
and
\begin{eqnarray}
\label{assumption-chi-2}
|\chi(y)| \leq C_- |\eta|, \quad \quad \quad \quad \quad y \in (-\infty,0).
\end{eqnarray}
Additionally, we assume the asymptotic behavior
\begin{equation}
\label{assumption-chi-3}
  \chi(y) \underset{y \to -\infty}{\sim} - \eps^{1/3} \eta.
\end{equation}
Under these assumptions on $\chi$, we consider the scalar equation
\begin{equation}
4 \nu''(y) + y \nu(y) - \nu^3(y) = \eps^{2/3} \left( 4 y \nu''(y) + 2 \nu'(y) +
\chi^2(y) \nu(y) \right), \quad y \in (-\infty,\eps^{-2/3}).  \label{P2limit-scalar}
\end{equation}
The Hastings--McLeod solution $\nu_0$ solves (\ref{P2limit-scalar})
for $\eps = 0$. For $\eps>0$ small, we are looking for a solution
$\nu$ to the scalar equation (\ref{P2limit-scalar}) near $\nu_0$. Thus,
using the decomposition $\nu = \nu_0 + R$, we rewrite (\ref{P2limit-scalar}) as
\begin{equation}
\label{persistence-problem}
M_{\eps} R = H_{\eps} + N_{\eps}(R),
\end{equation}
where the linear operator $M_{\eps}$, the source term $H_{\eps}$ and the nonlinear
function $N_{\eps}(R)$ are given by
\begin{equation}
\label{operator-L}
M_{\eps} = -4(1-\eps^{2/3}y)^{1/2} \partial_y (1-\eps^{2/3}y)^{1/2} \partial_y +W_0(y) ,
\end{equation}
\begin{equation}
\label{residual}
H_{\eps} := -\eps^{2/3} \left( 4 y \nu_0''(y) + 2 \nu_0'(y) + \chi^2(y) \nu_0(y) \right)
\end{equation}
and
\begin{equation}
\label{nonlinear}
N_{\eps}(R) := -3 \nu_0 R^2 - R^3 - \eps^{2/3} \chi^2 R,
\end{equation}
with $W_0(y) := 3 \nu_0^2(y) - y$. By Lemma 2.2 in \cite{GP},
there is a positive constant $W_{\rm min}$ such that
\begin{equation}
W_0(y) \geq W_{\rm min} \quad \mbox{\rm for all} \quad y \in \mathbb{R}.
\label{W-potential}
\end{equation}
Let us define the Hilbert spaces $L^2_{\eps}$ and $H^1_{\eps}$ as the sets of functions in
$L^1_{\rm loc}(-\infty,\eps^{2/3})$ with finite squared norms
\begin{eqnarray*}
\| u \|_{L^2_{\eps}}^2 & := & \int_{-\infty}^{\eps^{-2/3}} (1 - \eps^{2/3} y)^{-1/2} |u(y)|^2 dy, \\
\| u \|_{H^1_{\eps}}^2 & := & \int_{-\infty}^{\eps^{-2/3}} (1 - \eps^{2/3} y)^{-1/2}
\left[ 4 (1 - \eps^{2/3} y) |u'(y)|^2 + W_0(y) |u(y)|^2 \right] dy.
\end{eqnarray*}
By Lemma 2.3 in \cite{GP}, $M_\eps$ is defined as a self-adjoint
unbounded invertible operator on $L_\eps^2$ and for $\eps > 0$ small enough,
the inverse operator satisfies the $\eps$-independent bound
\begin{equation}
\label{bound-on-inverse-2}
\forall f \in L^2_\eps, \quad \| M_{\eps}^{-1} f \|_{H^1_{\eps}} \leq  W_{\rm min}^{-1/2} \| f \|_{L^2_{\eps}}.
\end{equation}
By the implicit function theorem arguments, we obtain the following result.

\begin{theo}
Let $\nu_0$ be the Hastings--McLeod solution of the Painlev\'e-II
equation, defined by (\ref{P2-eq}) and (\ref{P-bc}).
Let $\chi \in L^{\infty}(-\infty,\eps^{-2/3})$ satisfy
(\ref{assumption-chi-1})--(\ref{assumption-chi-3}).
For any $q > \frac{5}{6}$, there exist $\eps_q > 0$, $\eta_q > 0$, and $C_q > 0$
such that for every $\eps \in (0,\eps_q)$ and $|\eta| < \eta_q \eps^{q}$,
there exists a unique solution $R \in H^1_{\eps}$ of equation (\ref{persistence-problem})
such that
\begin{equation}
\label{upper-bound-nu}
\| \nu - \nu_0 \|_{H^1_{\eps}} \leq
C_q \left\{ \begin{array}{l} \eps^{2q - 4/3} \left| \log(\eps) \right|^{1/2}, \quad {\rm if} \; q \leq 1, \\
\eps^{2/3}, \quad  \quad  \quad  \quad  \quad  \quad  \quad {\rm if} \; q > 1. \end{array} \right.
\end{equation}
If $\nu := \nu_0 + R$, then $\nu(y)  > 0$ for all $y \in (-\infty,\eps^{-2/3})$ and
there is $\gamma > 0$ such that
\begin{equation}
\label{decay-nu}
\nu(y) \underset{y \to -\infty}{\sim} \gamma  |y|^{\frac{1-\eps - \eps^2 \eta^2}{4\eps}} e^{-\frac{|y|}{2 \eps^{1/3}}}.
\end{equation}
Furthermore, if $\nu_{1,2}$ correspond to $\chi_{1,2}$, then
there exists an $\eps$-independent positive constant $C$ such that
\begin{equation}
\label{Lipschitz-nu}
\| \nu_1 - \nu_2 \|_{H^1_{\eps}} \leq C \eps^{2/3} \| (\chi^2_1 - \chi_2^2) \nu_1 \|_{L^2_{\eps}}.
\end{equation}
\label{theorem-map-chi-to-nu}
\end{theo}
The proof of this theorem is divided into three subsections.

\subsection{Nonlinear and residual terms $N_{\eps}(R)$ and $H_{\eps}$}

First, we note the following embedding property.

\begin{lem}
There exists  $C>0$ such that if $\eps>0$ is small enough and if $u \in H^1_{\eps}$, then
$u \in \mathcal{C}^0(-\infty,\eps^{-2/3})$ satisfies
\begin{equation}
\label{embedding-1}
\| u \|_{L^{\infty}(-\infty,\eps^{-2/3})} \leq C \| u \|_{H^1_{\eps}}.
\end{equation}
\label{lemma-embeddings}
\end{lem}

\begin{Proof}
We introduce the map $T_\eps$ defined for $u\in L^2_\eps$ by
$$
(T_\eps u)(z)=u(\eps^{-2/3}-\eps^{2/3}z^2),\quad z\in\R.
$$
In \cite{GP}, we showed that $T_\eps$ is an isometry between $L^2_\eps$ and
the space $L^2_{\rm even}(\R)$ of even squared-integrable functions on $\R$.
Also, $T_{\eps}$ induces an isometry between $H^1_\eps$ and
\begin{eqnarray*}
H^1_w &=&\left\{f\in L^2_{\rm even}(\R) : \quad
\| f \|_{H^1_w}^2  :=  \int_{\mathbb{R}} \left( |f'(z)|^2 +  W_0(\eps^{-2/3}-\eps^{2/3} z^2) |f(z)|^2 \right) dz<\infty\right\} .
\end{eqnarray*}
As a result, Sobolev embedding implies for every $u\in H^1_\eps$ that
$$
\| u \|_{L^{\infty}(-\infty,\eps^{2/3})} = \| T_\eps u \|_{L^{\infty}(\mathbb{R})} \lesssim \| T_\eps u \|_{H^1(\R)} \lesssim \| T_\eps u \|_{H_w^1}=  \| u \|_{H^1_{\eps}},
$$
which yields (\ref{embedding-1}).
\end{Proof}

Next, we write $N_{\eps}(R) = N_0(R) + \Delta N_{\eps}(R)$, where
$$
N_0(R) = -3 \nu_0 R^2 - R^3, \quad
\Delta N_{\eps}(R) = -\eps^{2/3} \chi^2 R,
$$
with a given $\chi \in L^{\infty}(-\infty,\eps^{-2/3})$. We estimate the nonlinear terms in the following lemma.

\begin{lem}\label{lemma-N0-DN}
There exists $C>0$ such that for every $R \in H^1_{\eps}$, we have
$N_0(R)\in L^2_{\eps}$ with
\begin{equation}
\| N_0(R) \|_{L^2_{\eps}} \leq C \left( \| \nu_0 \|_{L^{\infty}(-\infty,\eps^{-2/3})} \| R \|^2_{H^1_{\eps}}
+ \| R \|^3_{H^1_{\eps}} \right).
\label{bound-nonlinear-term}
\end{equation}
There exists $C>0$ such that for $\eps>0$ small enough,
for $\chi \in L^{\infty}(-\infty,\eps^{-2/3})$ satisfying (\ref{assumption-chi-1}) and (\ref{assumption-chi-2}),
and for every $R \in H^1_{\eps}$, we have
\begin{equation}
\| \Delta N_{\eps}(R) \|_{L^2_{\eps}} \leq C \eps^{-2/3} \eta^2 \| R \|_{H^1_{\eps}}.
\label{bound-nonlinear-term-additional}
\end{equation}
\end{lem}

\begin{Proof}
By Sobolev embedding of $H^1(\mathbb{R})$ into $L^p(\mathbb{R})$ for any $p \geq 2$, for every $R\in H^1_\eps$ and $p=2,3$,
we have
$$
\| R^p \|_{L^2_{\eps}} = \| T_\eps(R^p) \|_{L^2}=\| T_\eps R \|_{L^{2p}}^p
\lesssim \|T_\eps R   \|^p_{H^1}\lesssim \|T_\eps R   \|^p_{H^1_w} =  \| R\|^p_{H^1_\eps},
$$
which yields bound (\ref{bound-nonlinear-term}).

Similarly, from (\ref{assumption-chi-1}) and (\ref{assumption-chi-2}), we have
\begin{equation*}
\| \Delta N_{\eps}(R) \|_{L^2_{\eps}} \leq \eps^{2/3} \| \chi \|^2_{L^{\infty}(-\infty,\eps^{-2/3})} \| R \|_{L^2_{\eps}}
\leq C \eps^{-2/3} \eta^2 \| R \|_{L^2_{\eps}},
\end{equation*}
which yields bound (\ref{bound-nonlinear-term-additional}).
\end{Proof}

\begin{rem}
\label{remark-constraints-1}
Since $\| \nu_0 \|_{L^{\infty}(-\infty,\eps^{-2/3})} = \mathcal{O}(\eps^{-1/3})$,
bound (\ref{bound-nonlinear-term}) implies that, for fixed $C_0 > 0$ and
$\alpha>1/3$, $N_0$ maps the ball of radius $C_0 \eps^\alpha$ centered at the
origin in $H^1_\eps$ into itself, provided $\eps$ is small enough. Moreover,
estimating $N_0(R_1)-N_0(R_2)$ similarly, one can show that if $\eps$ is
small enough, $N_0$ induces a contraction on these balls.
\end{rem}

\begin{rem}
\label{remark-constraints-2}
The term $\eps^{2/3} \chi^2$ in $\Delta N_{\eps}(R)$ is a small bounded perturbation to the linear operator $M_{\eps}$
if $\eta = \mathcal{O}(\eps^{q})$ with $q > \frac{1}{3}$, which is satisfied if $q > \frac{5}{6}$.
\end{rem}

Finally, we write $H_{\eps} = H_0 + \Delta H_{\eps}$, where
\begin{eqnarray*}
H_0 = -\eps^{2/3} \left( 4 y \nu_0''(y) + 2 \nu_0'(y) \right), \quad
\Delta H_{\eps} = -\eps^{2/3} \chi^2 \nu_0,
\end{eqnarray*}
with a given $\chi \in L^{\infty}(-\infty,\eps^{-2/3})$.  We estimate the residual terms in the following lemma.

\begin{lem}\label{lemma-N0-DN}
There exists $C>0$ such that
\begin{equation}
\label{bound-before-last}
\| H_0 \|_{L^2_{\eps}} \leq C \eps^{2/3},
\end{equation}
There exists $C>0$ such that for $\eps>0$ small enough and for $\chi \in L^{\infty}(-\infty,\eps^{-2/3})$
satisfying (\ref{assumption-chi-1}) and (\ref{assumption-chi-2}), we have
\begin{equation}
\label{bound-last}
\| \Delta H_{\eps} \|_{L^2_{\eps}} \leq
C  \eps^{-4/3} \eta^2 |\log(\eps)|^{1/2}.
\end{equation}
\end{lem}

\begin{Proof}
The first term $H_0$ was analyzed in \cite{GP}. The bound (\ref{bound-before-last}) holds
because $4 y \nu_0''(y) + 2 \nu_0'(y) \underset{y \to +\infty}{=} \mathcal{O}(y^{-7/2})$,
whereas this function decays even faster as $y\to -\infty$.

The second term is analyzed with the following auxiliary result,
\begin{eqnarray}
\int_0^{\eps^{-2/3}} \frac{d y}{(1 - \eps^{2/3} y)^{1/2} (1 + y)} &=&
\int_0^{\eps^{-2/3}} \frac{d y}{1 + y}+\int_0^{\eps^{-2/3}} \frac{1-(1 - \eps^{2/3} y)^{1/2} }{(1 - \eps^{2/3} y)^{1/2} (1 + y)}dy\nonumber\\
&=&\log(1+\eps^{-2/3})+\int_0^1\frac{1-(1-t)^{1/2}}{(1-t)^{1/2}(\eps^{2/3}+t)}dt\nonumber\\
&=&-\frac{2}{3}\log(\eps) +\mathcal{O}(1),
\label{bounds-on-integral}
\end{eqnarray}
where for the last equality, we have used Lebesgue's theorem, which is possible since for every $t\in(0,1)$ and $\eps>0$,
$$
\frac{1-(1-t)^{1/2}}{(1-t)^{1/2}(\eps^{2/3}+t)}\leq g(t) :=
\left\{\begin{array}{ll}\frac{t}{(1-t)^{1/2}t} \sup_{s\in[0,t]} \frac{1}{2\sqrt{1-s}} = \frac{1}{2(1-t)}& {\rm if\ }t\in (0,1/2)\\
\frac{2}{(1-t)^{1/2}} & {\rm if\ }t\in (1/2,1)\end{array}\right.
$$
and $g\in L^1(0,1)$. Hence, from (\ref{assumption-chi-1}) and (\ref{assumption-chi-2}), we have
$$
\| \Delta H_{\eps} \|_{L^2_{\eps}} \leq \eps^{2/3} \| \chi \|_{L^{\infty}(-\infty,\eps^{-2/3})}^2 \| \nu_0 \|_{L^2_{\eps}}
\leq C  \eps^{-2/3} \eta^2 \| \nu_0 \|_{L^2_{\eps}}
\leq C  \eps^{-4/3} \eta^2 \| (1+|y|)^{-1} \nu_0 \|_{L^2_{\eps}},
$$
where $(1+|y|)^{-1} \nu_0 = \mathcal{O}(y^{-1/2})$ as $y \to \infty$ and decays fast as $y \to -\infty$.
By using (\ref{bounds-on-integral}), this bound yields (\ref{bound-last}).
\end{Proof}

\subsection{Existence and properties of $R \in H^1_{\eps}$}

For $\eps > 0$ small enough, let $\chi \in L^{\infty}(-\infty,\eps^{-2/3})$ satisfy
(\ref{assumption-chi-1}), (\ref{assumption-chi-2}), and (\ref{assumption-chi-3}). Then,
we prove the existence of a unique solution
$R \in H^1_{\eps}$ of equation (\ref{persistence-problem}) satisfying
(\ref{upper-bound-nu}) provided that $\eta = \mathcal{O}(\eps^q)$
as $\eps \to 0$ for any $q > \frac{5}{6}$.

The existence of $R$  follows from a fixed-point argument in $\mathcal{B}_\eps$,
where $\mathcal{B}_\eps$ denotes the ball of $H^1_\eps$ centered at the origin, with radius
\begin{equation}
\label{radius-ball}
\rho_{\eps} := C_q \left\{ \begin{array}{l} \eps^{2q - 4/3} \left| \log(\eps) \right|^{1/2}, \quad {\rm if} \; q \leq 1, \\
\eps^{2/3}, \quad  \quad  \quad  \quad  \quad  \quad  \quad {\rm if} \; q > 1, \end{array} \right.
\end{equation}
for some $C_q > 0$. Indeed, inverting $M_{\eps}$, we rewrite (\ref{persistence-problem}) as
the fixed point equation
\begin{eqnarray}\label{recursion}
R = \Phi(R), \quad \Phi(R) := M_{\eps}^{-1} \left( H_0 + \Delta H_{\eps} + N_0(R) + \Delta N_{\eps}(R) \right),
\end{eqnarray}
By bounds (\ref{bound-on-inverse-2}), (\ref{bound-nonlinear-term}), (\ref{bound-nonlinear-term-additional}),
(\ref{bound-before-last}), and (\ref{bound-last}), we obtain
$$
\| \Phi(R) \|_{H^1_{\eps}} \leq C \left( \eps^{2/3} + \eps^{2q-4/3} |\log(\eps)|^{1/2} + \eps^{-1/3} \rho_{\eps}^2
+ \rho_{\eps}^3 + \eps^{2q-2/3} \rho_{\eps}\right).
$$
If $q > 1$, then $2q - \frac{4}{3} > \frac{2}{3}$ and for $\eps > 0$ small enough,
$\eps^{2/3} \gg \eps^{2q - 4/3} |\log(\eps) |^{1/2}$, thus
the operator $\Phi$ maps the ball $\mathcal{B}_\eps$ to itself.
If $q \leq 1$ but $q > \frac{5}{6}$, then $\eps^{2/3} \ll \eps^{2q - 4/3} |\log(\eps) |^{1/2}$ and
the operator $\Phi$ maps the ball $\mathcal{B}_\eps$ to itself.
Similarly, one can show that $\Phi$ is a contraction on the ball $\mathcal{B}_{\eps}$,
see Remarks \ref{remark-constraints-1} and \ref{remark-constraints-2}.

Next, we set $\nu := \nu_0 + R$ and prove positivity of $\nu(y)$ for all $y \in (-\infty,\eps^{-2/3})$
and decay of $\nu(y)$ as $y \to -\infty$, according to the asymptotic behavior (\ref{decay-nu}).

By Sobolev embedding (\ref{embedding-1}), we have $\| R \|_{L^{\infty}(-\infty,\eps^{-2/3})} \leq C \rho_{\eps}$,
where $\rho_{\eps} \to 0$ as $\eps \to 0$. Since $\nu_0$ is increasing and $\nu_0(0) > 0$, then
$\nu(y) > 0$ for all $y \in (0,\eps^{-2/3})$ and $\eps > 0$ small enough.
Additionally, we know that $R \in \mathcal{C}^0(-\infty,\eps^{-2/3})$ and $R(y) \to 0$
as $y \to -\infty$. By bootstrapping arguments, we obtain a higher regularity of
$\nu \in \mathcal{C}^2(-\infty,\eps^{-2/3})$, and hence $\nu \in \mathcal{C}^{\infty}(-\infty,\eps^{-2/3})$.

Next, coming back to the variable $x$ in the transformation (\ref{asymptotic-transformation}),
we can see that $\varphi(x) = \eps^{1/3} \nu(y)$ satisfies the first equation of system (\ref{stationaryGP}), whereas
since $\chi$ satisfies (\ref{assumption-chi-3}), then
$$
\xi \underset{|x| \to \infty}{\sim} - \eps \eta.
$$
By the same method as in the proof of Lemma \ref{lemma-ode},
we obtain positivity of $\nu(y)$ for $y \in (-\infty,0)$ and the decay of
$\nu(y) \to 0$ as $y \to -\infty$.
The decay behavior (\ref{decay-nu}) follows from
the decay behavior (\ref{constraint-1}) by the change of variables (\ref{asymptotic-transformation}).

\subsection{Lipschitz continuity of the map $\chi \mapsto \nu$}

We prove the bound (\ref{Lipschitz-nu}) and hence complete the proof of Theorem \ref{theorem-map-chi-to-nu}.
First, we write equation (\ref{P2limit-scalar}) for $\nu_1$ and $\nu_2$ related to $\chi_1$ and $\chi_2$.
Taking the difference and denoting $\delta \nu := \nu_1 - \nu_2$, we obtain
\begin{equation}
\label{P2limit-scalar-difference}
\left( M_{\eps} + \Delta W_1 + \Delta W_2 \right) \delta \nu = - \eps^{2/3} (\chi_1^2 - \chi_2^2) \nu_1,
\end{equation}
where $\Delta W_1 := \eps^{2/3} \chi_2^2$ and $\Delta W_2 := \nu_1^2 + \nu_1 \nu_2 + \nu_2^2 - 3 \nu_0^2$.
By the assumptions (\ref{assumption-chi-1}) and (\ref{assumption-chi-2}), there is
an $\eps$-independent positive constant $C$ such that
$$
\| \Delta W_1 \|_{L^{\infty}(-\infty,\eps^{-2/3})} \leq C \eta^2 \eps^{-2/3},
$$
which shows that, if $\eta = \mathcal{O}(\eps^{q})$ with $q > \frac{1}{3}$ as $\eps \to 0$,
then $\Delta W_1$ is a small bounded perturbation to the positive potential $W_0$ in $M_{\eps}$.
On the other hand, denoting $\nu_{1,2} = \nu_0 + R_{1,2}$, we have
$$
\Delta W_2 = 3 \nu_0 (R_1 + R_2) + R_1^2 + R_1 R_2 + R_2^2.
$$
Since both $R_1$ and $R_2$ belongs to $\mathcal{B}_\eps$ with radius (\ref{radius-ball}),
there is another $\eps$-independent positive constant $C$ such that
$$
\| \Delta W_2 \|_{L^{\infty}(-\infty,\eps^{-2/3})} \leq C \rho_{\eps} \eps^{-1/3}.
$$
If $q > \frac{5}{6}$, then $\Delta W_2$ is another small bounded perturbation
to the positive potential $W_0$ in $M_{\eps}$. Hence
$\left( M_{\eps} + \Delta W_1 + \Delta W_2 \right)$
is an invertible operator with an $\eps$-independent bound on its inverse
from $L^2_{\eps}$ to $H^1_{\eps}$. Therefore, we obtain
from equation (\ref{P2limit-scalar-difference}) that
$$
\| \delta \nu \|_{H^1_{\eps}} = \eps^{2/3} \| \left( M_{\eps} + \Delta W_1 + \Delta W_2 \right)^{-1} (\chi^2_1 - \chi_2^2) \nu_1 \|_{H^1_{\eps}}
\leq C \eps^{2/3} \| (\chi^2_1 - \chi_2^2) \nu_1 \|_{L^2_{\eps}},
$$
which yields the bound (\ref{Lipschitz-nu}).

\section{Mapping $\nu \to \chi$}

Here we consider the integral formula (\ref{chi}) for a family of positive
functions $\nu \in L^2(-\infty,\eps^{-2/3})\cap \mathcal{C}^0(-\infty,\eps^{-2/3})$,
which depends on $\eps$. This integral formula defines a solution of the second equation in the system
(\ref{P2limit}). We assume that there are constants $C_+ \geq 1$ and $C_- > 0$ such that
for every $\eps > 0$ small enough, the function $\nu$ satisfies
\begin{eqnarray}
\label{assumption-nu-1}
\frac{1}{C_+} y \leq \nu^2(y) \leq C_+ (1+y), \quad y \in (0,\eps^{-2/3}),
\end{eqnarray}
and
\begin{eqnarray}
\label{assumption-nu-2}
\frac{\int_{-\infty}^y \nu^2(s) ds}{\nu^2(y)} \leq C_-, \quad y \in (-\infty,0).
\end{eqnarray}
In addition, we assume that there is $\gamma > 0$ such that $\nu$ satisfies the asymptotic decay
\begin{equation}
\label{assumption-nu-3}
\nu(y) \underset{y \to -\infty}{\sim} \gamma |y|^{\frac{1-\eps - \eps^2 \eta^2}{4\eps}} 
e^{-\frac{|y|}{2 \eps^{1/3}}}.
\end{equation}

We shall study the mapping $\nu \to \chi$, defined
on some neighborhood of $\nu_0$ in a suitable
space such that (\ref{chi}) provides a bounded function $\chi$.
First, we obtain the following elementary result.

\begin{lem}
Let $\nu \in L^2(-\infty,\eps^{-2/3})\cap \mathcal{C}^0(-\infty,\eps^{-2/3})$ satisfy
(\ref{assumption-nu-1})--(\ref{assumption-nu-3}).
Then, $\chi \in L^{\infty}(-\infty,\eps^{-2/3})$ is well-defined by the integral formula (\ref{chi}) and
satisfies properties (\ref{assumption-chi-1})--(\ref{assumption-chi-3}).
\label{theorem-map-nu-to-chi}
\end{lem}

\begin{Proof}
Thanks to assumptions (\ref{assumption-nu-1}) and (\ref{assumption-nu-2}),
we obtain
$$
\frac{(1+y)^2-1}{2 C_+} \leq
\| \nu \|_{L^2(-\infty,y)}^2 \leq
C_+ \left(C_- + \frac{(1+y)^2-1}{2}\right) , \quad  y \in [0,\eps^{-2/3}].
$$
Hence, the lower and upper bounds on $|\chi(y)|$ in (\ref{assumption-chi-1})
follow from the lower and upper bounds on $\nu^2(y)$ in (\ref{assumption-nu-1}).
Bound (\ref{assumption-chi-2}) follows from the definition (\ref{chi}) and bound (\ref{assumption-nu-2}).
Finally, the asymptotic decay (\ref{assumption-nu-3}) gives
$$
\frac{\int_{-\infty}^y \nu^2(s) ds}{\nu^2(y)} \underset{y \to -\infty}{\sim} \eps^{1/3},
$$
which is equivalent to the property (\ref{assumption-chi-3}).
\end{Proof}

The following result gives a Lipschitz continuity property of the mapping $\nu \to \chi$
in a neighborhood of $\nu_0$.

\begin{lem}
\label{lemma-mapping}
Let $\chi_{1,2}$ be defined by (\ref{chi}) for $\nu_{1,2} \in L^2(-\infty,\eps^{-2/3})\cap \mathcal{C}^0(-\infty,\eps^{-2/3})$,
where $\nu_{1,2}$ satisfy (\ref{assumption-nu-1})--(\ref{assumption-nu-2}) and are close to $\nu_0$ so that
for a positive $\delta$, they satisfy the bound
$$
\| \nu_{1,2} - \nu_0 \|_{L^2(-\infty,\eps^{-2/3})} + \| \nu_{1,2} - \nu_0 \|_{L^{\infty}(-\infty,\eps^{-2/3})} \leq \delta.
$$
For $\eps > 0$ sufficiently small, there is an $\eps$-independent positive constant $C$ such that
\begin{eqnarray}
\label{chi-full-bound-2}
\| \chi_1 - \chi_2 \|_{L^{\infty}(0,\eps^{-2/3})} \leq C |\eta| \left( \| \nu_1 - \nu_2 \|_{L^2(-\infty,\eps^{-2/3})}
+ \eps^{-1/3}\| \nu_1 - \nu_2 \|_{L^{\infty}(0,\eps^{-2/3})}  \right).
\end{eqnarray}
Furthermore, for any fixed $y_0 \in (-\infty,0)$ and $\eps > 0$ sufficiently small,
there is another $\eps$-independent positive constant $C(y_0)$ such that
\begin{equation}
\label{bound-Lipschitz-1}
\|\chi_1-\chi_2\|_{L^\infty(y_0,0)} \leq C(y_0) |\eta|
\left( \| \nu_{1} - \nu_2 \|_{L^2(-\infty,\eps^{-2/3})} + \| \nu_{1} - \nu_2 \|_{L^{\infty}(-\infty,\eps^{-2/3})} \right).
\end{equation}
\end{lem}

\begin{Proof}
We write
\begin{equation}
\label{bounds-technical-2}
\chi_1(y) - \chi_2(y) = -\frac{\eta}{\nu_1^2(y)} \int_{-\infty}^y (\nu_1^2(s) - \nu_2^2(s)) ds
- \frac{\chi_2(y)}{\nu_1^2(y)}  (\nu_1^2(y) - \nu_2^2(y)).
\end{equation}
By using the Cauchy--Schwarz and triangle inequalities, we obtain
\begin{eqnarray}
\nonumber
| \chi_1(y) - \chi_2(y) | & \leq &
\frac{|\chi_1(y)|}{\| \nu_1 \|_{L^2(-\infty,y)}}
\left(1 + \frac{\| \nu_2 \|_{L^2(-\infty,y)}}{\| \nu_1 \|_{L^2(-\infty,y)}}\right)
\| \nu_1 - \nu_2 \|_{L^2(-\infty,y)} \\
\label{bounds-technical-0}
& \phantom{t} & \phantom{texttext} +
\frac{|\chi_2(y)|}{\nu_1^2(y)} (\nu_1(y)+ \nu_2(y)) |\nu_1(y) -
\nu_2(y) |.
\end{eqnarray}
Since $\nu_0$ is increasing, we have $\nu_0(0) > 0$ and $\| \nu_0 \|_{L^2(-\infty,0)} > 0$.
As a result, for $\delta > 0$ sufficiently small, we have
$$
\underset{y\in[0,\eps^{-2/3}]}{\inf}\nu_{1,2}(y) > \frac{1}{2} \nu_0(0)
$$
and
$$
\| \nu_{1,2} \|_{L^2(-\infty,0)} > \frac{1}{2} \| \nu_{0} \|_{L^2(-\infty,0)}.
$$
Using moreover (\ref{assumption-chi-1}) and (\ref{assumption-nu-1})
in the bound (\ref{bounds-technical-0}) for $y \in (0,\eps^{-2/3})$, we obtain (\ref{chi-full-bound-2}).

Using (\ref{bounds-technical-0}) for $y \in (-\infty,0)$ and
(\ref{assumption-chi-2}), we obtain for $y\in (y_0,0)$
\begin{eqnarray*}
\|\chi_1 - \chi_2 \|_{L^{\infty}(y_0,0)} & \leq & C_- |\eta|
\left( 1 + \frac{\| \nu_2 \|_{L^2(-\infty,0)}}{\| \nu_1 \|_{L^2(-\infty,y_0)}} \right)
\frac{\| \nu_1 - \nu_2 \|_{L^2(-\infty,0)}}{\| \nu_1 \|_{L^2(-\infty,y_0)}} \\
& \phantom{t} & + C_- |\eta|
\left(1 + \frac{\sup_{y \in [y_0,0]} |\nu_2(y)|}{\inf_{y \in [y_0,0]} |\nu_1(y)|} \right)
\frac{\| \nu_1 - \nu_2\|_{L^{\infty}(y_0,0)}}{\inf_{y \in [y_0,0]} |\nu_1(y)|},
\end{eqnarray*}
Provided
$$
\delta < \frac{1}{2} \min\left[\|\nu_0\|_{L^2(-\infty,0)},\nu_0(y_0)\right],
$$
the triangle inequality $|\nu_1 | \geq |\nu_0| - |\nu_1 - \nu_0|$ yields
\begin{eqnarray}\label{2}
\|\chi_1 - \chi_2 \|_{L^{\infty}(y_0,0)} & \leq & C_- |\eta|
\left( 1 + \frac{\| \nu_0 \|_{L^2(-\infty,0)}+\delta}{\| \nu_0 \|_{L^2(-\infty,y_0)} - \delta} \right)
\frac{\| \nu_1 - \nu_2 \|_{L^2(-\infty,\eps^{-2/3})}}{\| \nu_0 \|_{L^2(-\infty,y_0)}-\delta} \nonumber\\
& \phantom{t} & + C_- |\eta|
\left(1 + \frac{\| \nu_0 \|_{L^{\infty}(-\infty,0)}+\delta}{\nu_0(y_0) - \delta} \right)
\frac{\| \nu_1 - \nu_2 \|_{L^{\infty}(-\infty,\eps^{-2/3})}}{\nu_0(y_0) - \delta},
\end{eqnarray}
which yields (\ref{bound-Lipschitz-1}).
\end{Proof}

\begin{rem}
Unless the Lipschitz continuity of Lemma \ref{lemma-mapping} is extended for $y_0 \to -\infty$
under the decay condition (\ref{assumption-nu-3}) with a good bound on the Lipschitz constant,
it is problematic to prove convergence of an iterative method for obtaining solutions of the system (\ref{P2limit})
by coupling the maps $\chi \to \nu$ and $\nu \to \chi$ together.
\end{rem}

\section{Solution in Conjecture \ref{theorem-2} via a numerical iterative method}

We shall develop an iterative numerical scheme to illustrate the validity
of the existence result stated in Conjecture \ref{theorem-2}.
The leading-order solution for the component $\nu$ of the problem (\ref{P2limit}) is the
Hastings--McLeod solution $\nu_0$ of the Painlev\'e-II equation defined by
(\ref{P2-eq})-(\ref{P-bc}). Let us define the zero iteration
for the component $\chi$ by
\begin{equation}
\label{chi-0}
\chi_0(y) = -\frac{\eta}{\nu_0^2(y)} \int_{-\infty}^y \nu_0^2(s) ds.
\end{equation}
Properties of this function are described by the following proposition.

\begin{prop}
$\chi_0 \in \mathcal{C}^\infty(\mathbb{R})$ and satisfies
the asymptotic behavior
\begin{equation}
\label{chi-bc}
\chi_0(y) = -\eta \left\{ \begin{array}{l}
\frac{1}{2} y + \frac{3}{2} y^{-2} + \mathcal{O}(y^{-5}) \;\;\; \mbox{as} \quad y\to +\infty \\
|y|^{-1/2} + \mathcal{O}(|y|^{-5/4}) \quad \mbox{as} \quad y\to-\infty \end{array} \right.
\end{equation}
\label{lemma-chi}
\end{prop}

\begin{Proof}
First, since $\nu_0(y) > 0$ for all $y \in \mathbb{R}$ \cite{HM,HM2}
and $\nu_0^2(y)$ decays fast as $y \to -\infty$, the integral formula (\ref{chi-0}) defines
$\chi_0(y)$ for every $y \in \mathbb{R}$. Moreover, since
$\nu_0 \in \mathcal{C}^{\infty}(\R)$, then $\chi_0 \in \mathcal{C}^{\infty}(\mathbb{R})$.
We shall now consider the asymptotic behavior of $\chi_0$ as $y \to \pm \infty$.

For $y \to -\infty$, we use the asymptotic behavior of $\nu_0$ given by (\ref{P-bc}) and integration by parts
to obtain
\begin{eqnarray*}
\int_{-\infty}^y \nu_0^2(s) ds & = &
\frac{1}{\pi} \int_{-\infty}^y |s|^{-1/2} e^{-\frac{2}{3} |s|^{3/2}} \left(1 + \mathcal{O}(|s|^{-3/4}) \right) ds \\
& \underset{y\to -\infty}{=} & \frac{1}{\pi |y|}  e^{-\frac{2}{3}|y|^{3/2}}
\left(1 + \mathcal{O}(|y|^{-3/4}) \right).
\end{eqnarray*}
Dividing this expression by $\nu_0^2(y)$ and using the asymptotic behavior (\ref{P-bc}) as $y \to -\infty$,
we obtain the second line of (\ref{chi-bc}).

For $y \to +\infty$, we use the asymptotic behavior (\ref{P-bc}) and write
\begin{eqnarray}
\label{large-y-behavior}
\int_{-\infty}^y \nu_0^2(s) ds & \underset{y\to +\infty}{=} & \frac{1}{2} y^2 + y^{-1} + \mathcal{O}(y^{-4}).
\end{eqnarray}
Dividing this expression by $\nu_0^2(y)$, we obtain the first line of (\ref{chi-bc}).
\end{Proof}

Replacing $\chi$ by $\chi_0$ in the scalar equation (\ref{P2limit-scalar}), we obtain the first
iteration $\nu_1$ from Theorem \ref{theorem-map-chi-to-nu}.
Because the asymptotic behavior (\ref{assumption-chi-3}) is replaced by
the asymptotic behavior given in the second line
of (\ref{chi-bc}), the asymptotic decay (\ref{decay-nu}) is modified as follows:
\begin{equation}
\label{decay-nu-mod}
\nu_1(y) \underset{y \to -\infty}{\sim} \gamma_1 |y|^{\frac{1-\eps}{4\eps}} e^{-\frac{|y|}{2 \eps^{1/3}}}.
\end{equation}
Nevertheless, Lemma \ref{theorem-map-nu-to-chi} is applied in spite
of the modification (\ref{decay-nu-mod}) to produce the first iterate
$\chi_1$ satisfying (\ref{assumption-chi-1})--(\ref{assumption-chi-3}).
Then, we compute the second iterates $\nu_2$ from Theorem \ref{theorem-map-chi-to-nu} and
$\chi_2$ from Lemma \ref{theorem-map-nu-to-chi}, and continue on this computational algorithm.

We will now implement this iterative scheme numerically to show that
the sequence $\{ (\nu_n,\chi_n) \}_{n \in \mathbb{N}}$ converges
to a solution of the coupled system (\ref{P2limit}).

\begin{figure}[h]
\begin{center}
\includegraphics[scale=0.33]{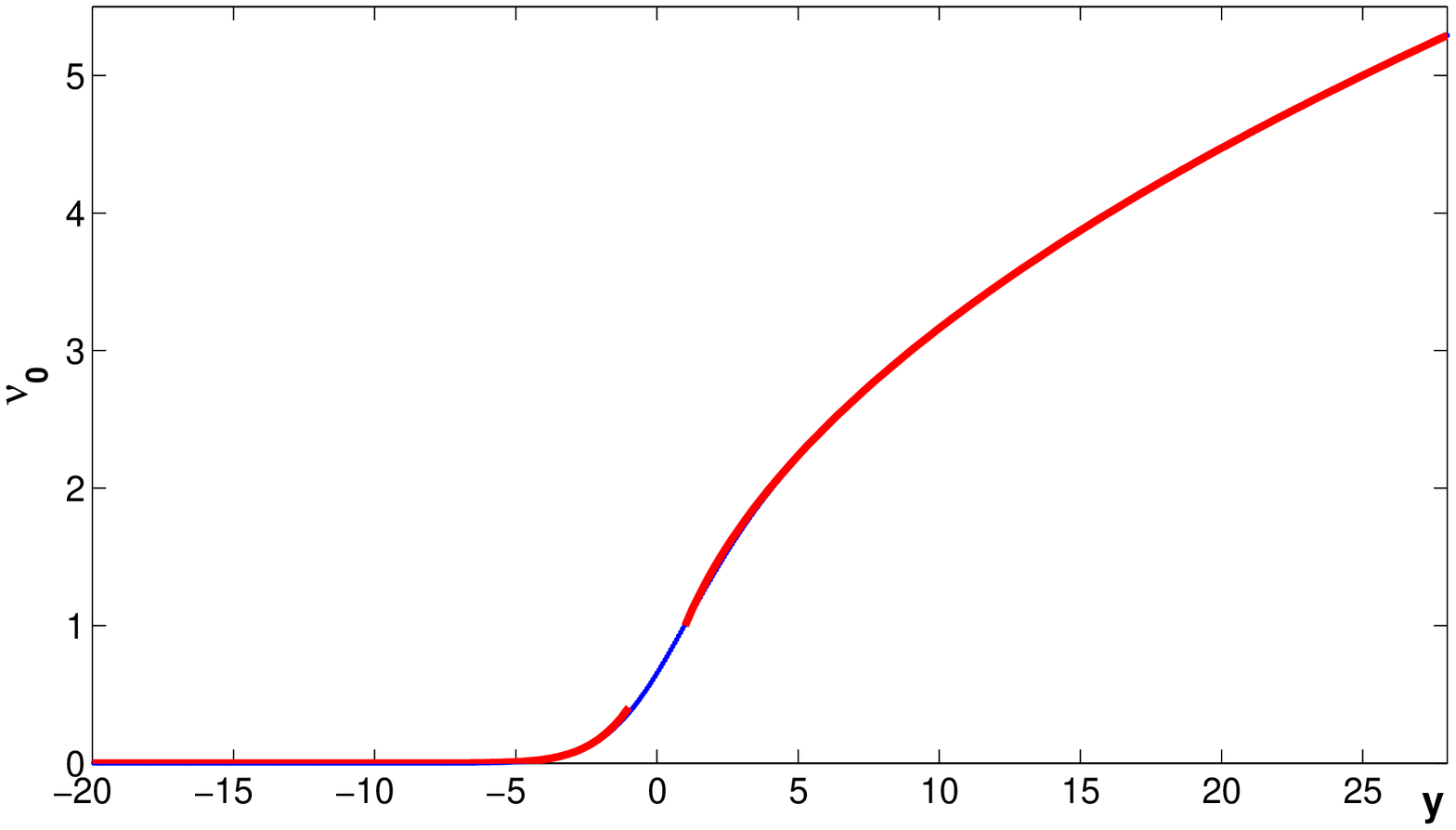}
\includegraphics[scale=0.32]{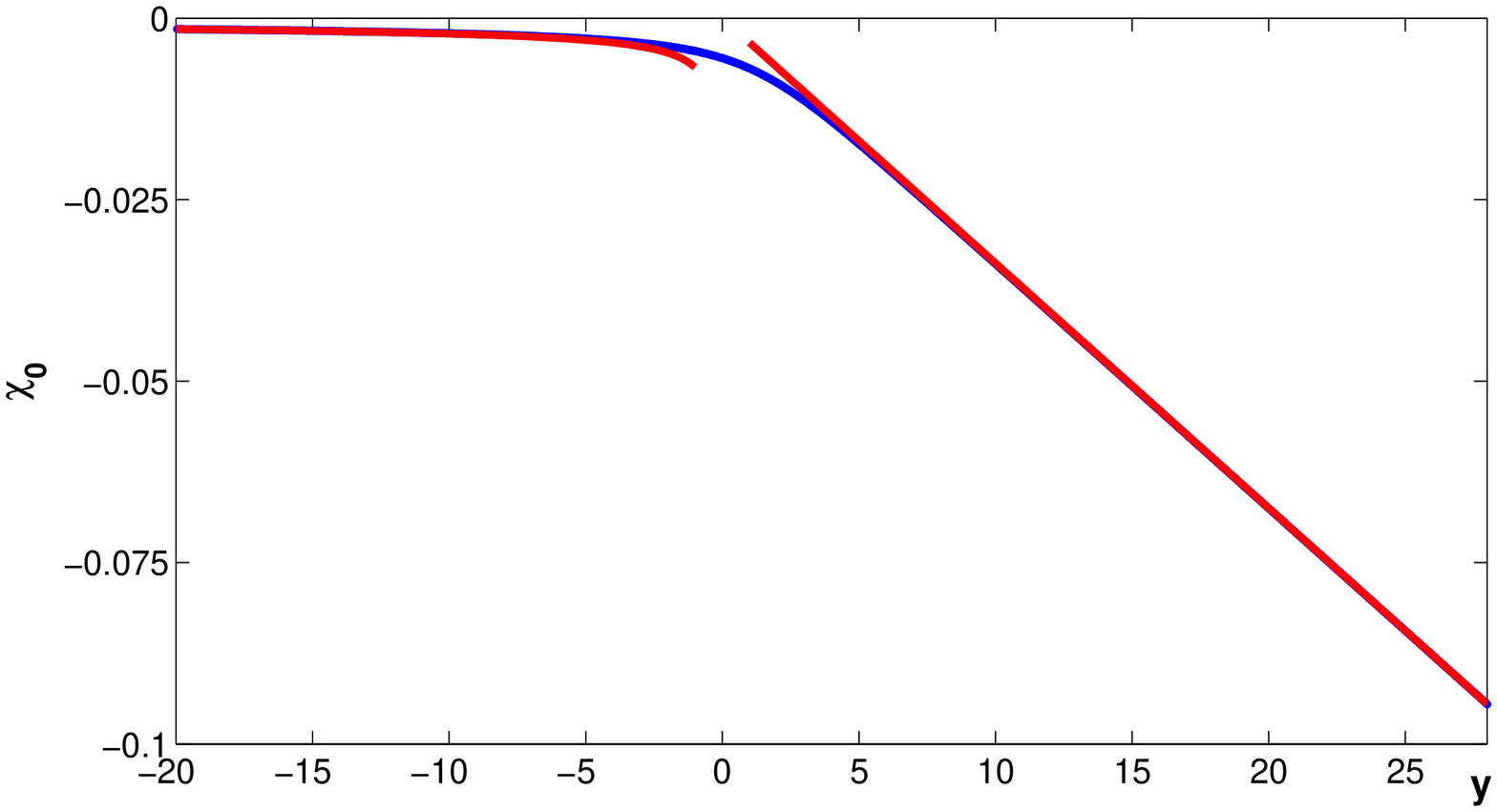}
\end{center}
\caption{\label{fig-Zeroprofile}
Components $\nu_0$ (left) and $\chi_0$ (right) for the numerical approximation
of the Hastings--McLeod solution of the Painlev\'e-II equation (\ref{P2-eq})
and the integral equation (\ref{chi-0}) for $\eta = \eps$ with $\eps = 0.0067$.}
\end{figure}

First, we approximate numerically the Hastings--McLeod solution
$\nu_0$ of the Painlev\'e-II equation (\ref{P2-eq}). We use the second-order
Heun's method supplemented with a shooting algorithm. The solution
is shown on the left panel of Figure \ref{fig-Zeroprofile}. The dashed lines showing
asymptotical expansions (\ref{P-bc}) for $y \geq 1$ and $y \leq -1$ are not
distinguished from the numerical approximations (dots).
We truncate the solution at $y_0 = 28$ and choose $\eps = y_0^{-3/2} = 0.0067$.
Then, we use the composite trapezoidal rule and approximate
the component $\chi_0$ from the integral equation (\ref{chi-0}) for $\eta = \eps$.
The solution is shown on the right panel of Figure \ref{fig-Zeroprofile}.
Again, the dashed lines show asymptotical expansions (\ref{chi-bc}) for $y \geq 1$ and $y \leq -1$.

Next, we use the iterative method to obtain the sequence $\{ (\nu_n,\chi_n) \}_{n \in \mathbb{N}}$
numerically. At each step, the numerical solution for $\nu_n$ is obtained from
the scalar equation (\ref{P2limit-scalar}) with $\chi = \chi_{n-1}$ by
the result of Theorem \ref{theorem-map-chi-to-nu}. Implemented numerically with a second-order difference method,
it takes just very few iterations to obtain a suitable approximation for $\nu_n$.
Then, the numerical solution for $\chi_n$ is obtained from the integral equation (\ref{chi})
with $\nu = \nu_n$ by applying the composite trapezoidal rule.
The iterations are terminated when the difference between two subsequent approximations
becomes smaller than $10^{-15}$. For the same value of $\eps = 0.0067$,
the numerical method converges in $3$ iterations
for $\eta = \eps$, in $6$ iterations for $\eta = \eps^{0.5}$, and in $11$ iterations
for $\eta = \eps^{0.25}$. No convergence of this method was found for $\eta = \eps^{0.15}$.

\begin{figure}[h]
\begin{center}
\includegraphics[scale=0.33]{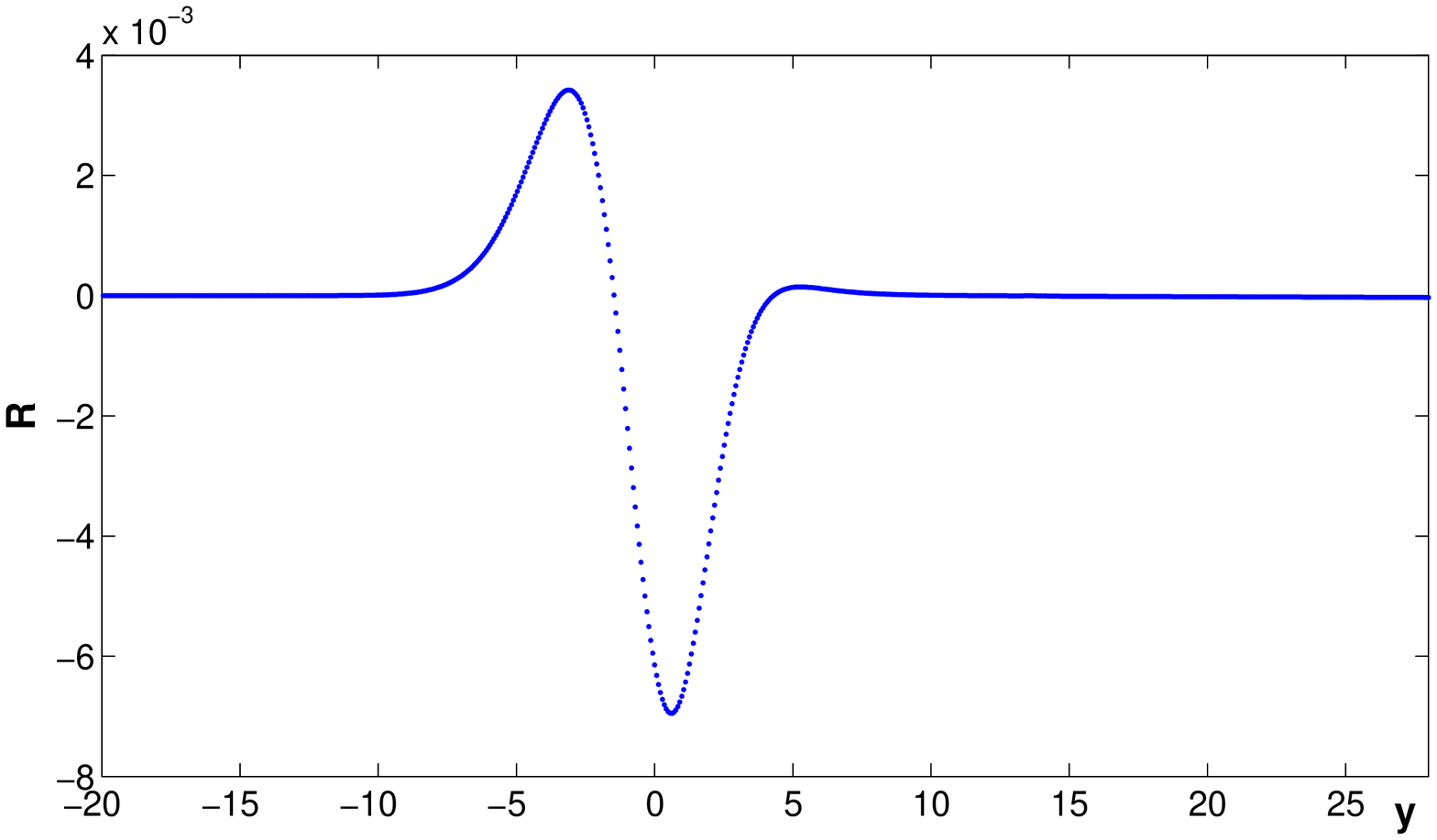}
\includegraphics[scale=0.32]{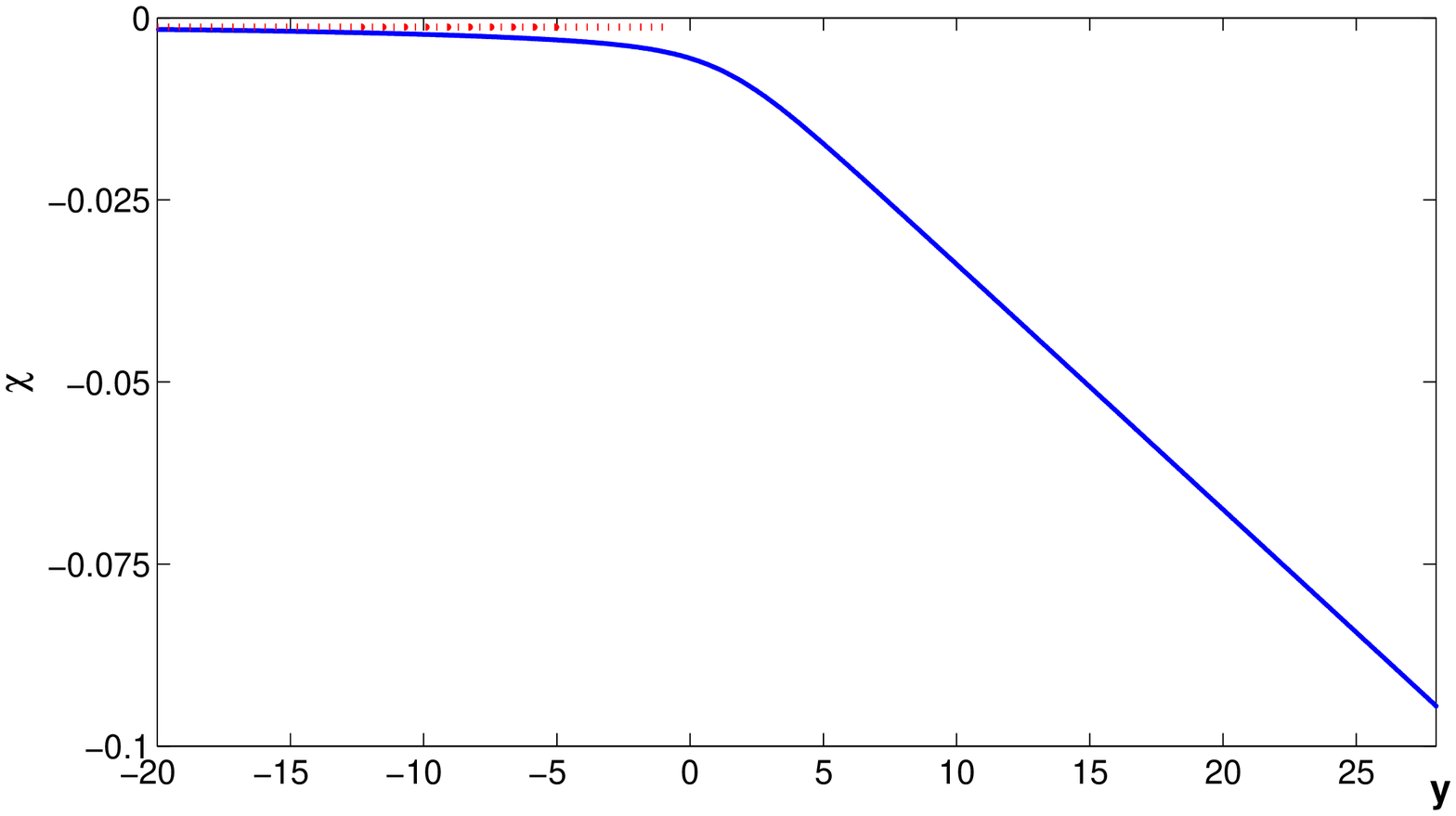} \\
\includegraphics[scale=0.33]{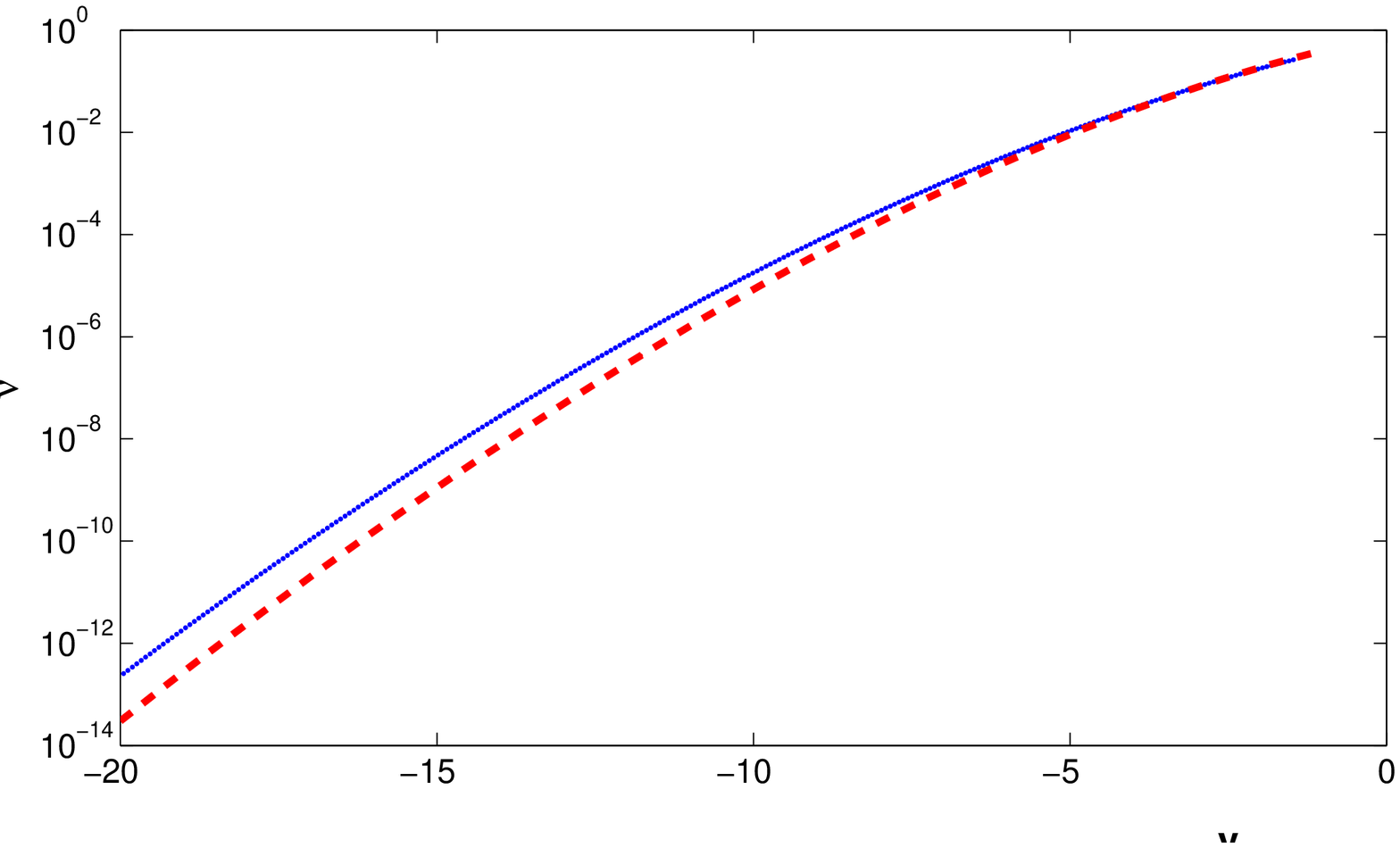}
\includegraphics[scale=0.32]{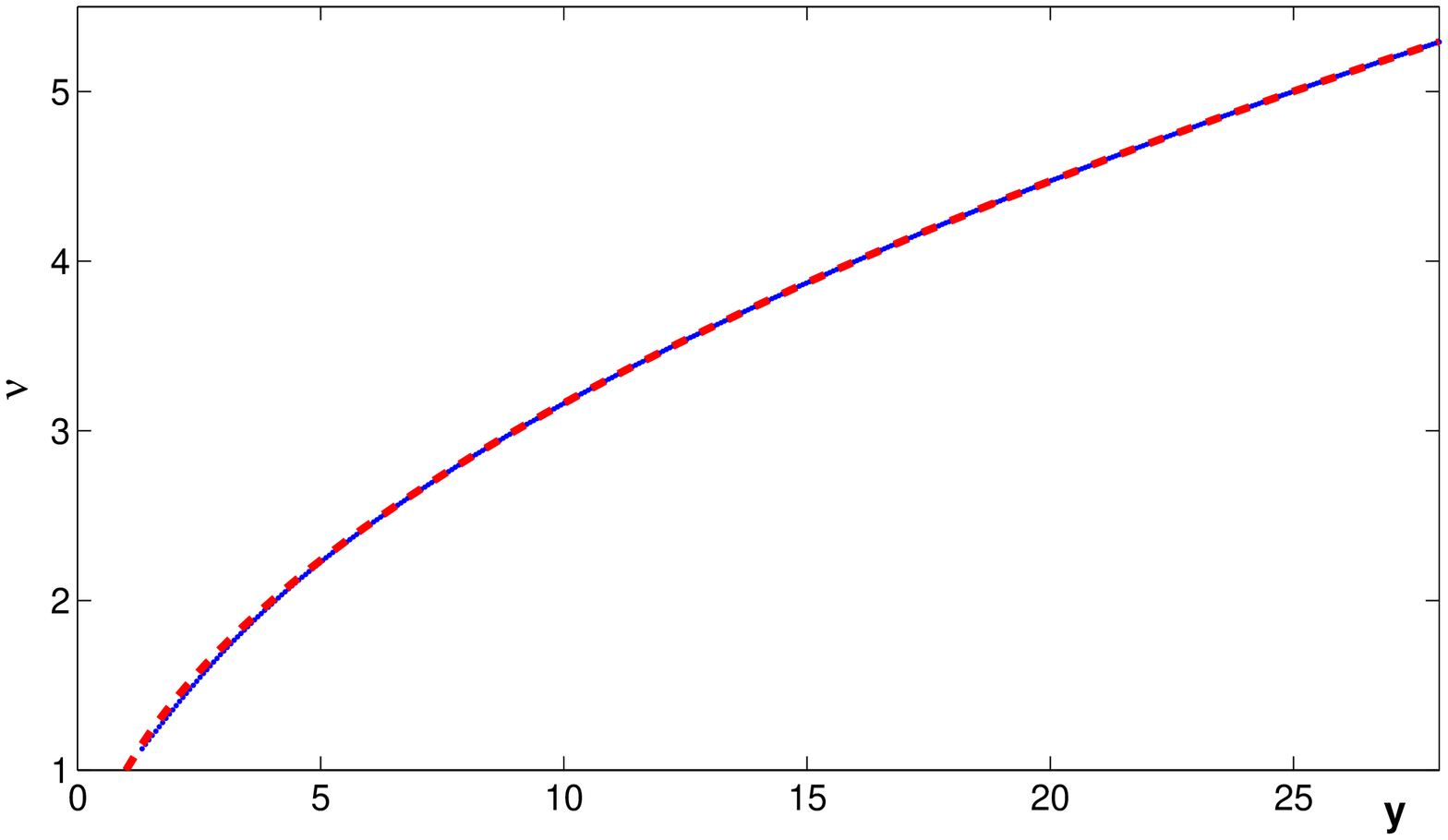}
\end{center}
\caption{\label{fig-Iterations}
Details of the numerical approximation of the
solution of the coupled system (\ref{P2limit})
for $\eta = \eps$ with $\eps = 0.0067$: component $R$ (top left panel),
component $\chi$ (top right panel), component $\nu$ (bottom panels) in comparison with
various asymptotic values shown by dashed lines.}
\end{figure}

Figure \ref{fig-Iterations} shows details of the numerical solution for $\eta = \eps$.
The top left panel shows the component $R := \nu - \nu_0$ of the final iterate of
the numerical solution. The top right panel shows the component $\chi$, where the dashed line
indicates the asymptotic value (\ref{assumption-chi-3}) for large negative $y$.
The bottom left panel shows the component $\nu$ (dots) in comparison with
the asymptotic decay behavior (\ref{P-bc}) of the Hastings--McLeod solution $\nu_0$ (dashed line).
It is clear from the semi-logarithmic scale that the component $\nu$ decays slower,
which agrees with the asymptotic behavior (\ref{decay-nu}). The bottom right panel shows
the component $\nu$ (dots) in comparison with the growth condition (\ref{P-bc})
of the Hastings--McLeod solution $\nu_0$ (dashed line). Because the values of
$R$ are small for $y = \mathcal{O}(\eps^{-2/3})$, the components $\nu$ and $\nu_0$
have similar growth rate. The situation changes when the value of $\eta$ is
larger, e.g. for $\eta = \eps^{0.25}$, when the values of $R$ become large
near the end $y = y_0$ of the computational interval. For such large values of $\eta$,
the bound (\ref{upper-bound-nu}) cannot be justified.

\section{Discussion}

Here we discuss the Thomas--Fermi limit of the ground state
in the stationary Gross--Pitaevskii equation with a more general $\PT$-symmetric potential:
\begin{equation}
\label{stationaryGP1-last}
\mu U(X) = \left( -\partial_X^2 + X^2 + 2 i \alpha X W(X) + |U(X)|^2\right) U(X),  \quad
X \in \R,
\end{equation}
where $W$ is an even, bounded, and decaying potential. In particular, we assume that
$$
X^n W(X) \in L^{\infty}(\mathbb{R}) \quad \mbox{\rm for all} \;\; n \in \mathbb{N}.
$$
Changing the variables $\mu = \eps^{-1}$, $x = \eps^{1/2} X$, and $U(X) = \eps^{-1/2} \varphi(x) e^{i \theta(x)}$,
and using the scaled variable $\xi(x) := \eps \theta'(x)$ and scaled
parameter $\alpha = \eps^{1/2} \eta$, we obtain the existence problem in the form,
\begin{equation}
\label{stationaryGP-last}
\left\{ \begin{array}{l}
\left( 1 - x^2 - \varphi^2(x) - \xi^2(x) \right) \varphi(x) = - \eps^2 \varphi''(x),  \\
\left( \varphi^2 \xi \right)'(x) = 2 \eta x W(\eps^{-1/2} x) \varphi^2(x), \end{array} \right. \quad x \in \R.
\end{equation}
The existence problem has now two scales $x$ and $X = \eps^{-1/2} x$ thanks to the bounded and
decaying potential $W$. As a result, the analysis of this existence problem at least
for finite and even large values of $\eta$ can be performed by a straightforward asymptotic method.

Solving the second equation of system (\ref{stationaryGP-last}) uniquely
from the condition $\lim\limits_{x \to \pm \infty} \varphi^2(x) \xi(x) = 0$, we obtain
an integral representation
\begin{equation}
\label{component-xi-last}
\xi(x) = \frac{2 \eta \eps}{\varphi^2(x)} \int_{-\infty}^{\eps^{-1/2} x} s W(s) \varphi^2(\eps^{1/2} s) ds.
\end{equation}
Assuming now that $\varphi \in \mathcal{C}^2$ near $x = 0$ and $X^3 W(X) \in L^1(\mathbb{R})$, we expand (\ref{component-xi-last})
into the asymptotic approximation,
\begin{equation}
\label{component-xi-approximation}
\xi(x) = 2 \eta \eps \left( \int_{-\infty}^{\eps^{-1/2} x} s W(s) ds \right) \frac{\varphi^2(0) + \mathcal{O}(\eps)}{\varphi^2(x)}.
\end{equation}
This asymptotic approximation shows that the phase-related component $\xi$ gives a contribution
to the $\PT$-symmetric ground state only if $\eta$ is as large as $\mathcal{O}(\eps^{-1})$ as $\eps \to 0$
and that this contribution is only affecting the ground state in the tiny region $|x| = \mathcal{O}(\eps^{1/2})$
around the origin as $\eps \to 0$. Therefore, the solution $\varphi$ of the existence problem (\ref{stationaryGP-last})
for $\eta = \mathcal{O}(\eps^{-1})$ is close to the solution of the existence problem (\ref{stationaryGP-last})
with $\eta = 0$ (which was justified in our previous work \cite{GP}), except for the values $|x| = \mathcal{O}(\eps^{1/2})$,
where the solution $\varphi$ is close to the modified Thomas--Fermi approximation
\begin{equation}
\label{TFlimit-last}
\varphi^2_{\rm TF}(x) = 1 - x^2 -  4 \eta^2 \eps^2 \left( \int_{-\infty}^{\eps^{-1/2} x} s W(s) ds \right)^2, \quad |x| \leq C \eps^{1/2},
\end{equation}
where $C$ is $\eps$-independent. From the requirement $\varphi_{\rm TF}^2(0) > 0$, we find the existence interval
$\eta \in (-\eta_0,\eta_0)$ of the $\PT$-symmetric ground state at the Thomas--Fermi limit, where
$$
\eta_0 := \frac{1}{2 \eps |\int_{-\infty}^0 s W(s) ds|}.
$$
Note that the breakdown of the ground state occurs at the origin $x = 0$, because
the absolute value of the integral $\int_{-\infty}^{\eps^{-1/2} x} s W(s) ds$
quickly drops when $x$ deviates from the origin. Therefore, we reiterate
the two facts mentioned in Remark \ref{remark-1}: the Thomas--Fermi radius $|x| = 1$ is
independent of the gain-loss parameter $\eta$ and the $\PT$-symmetric potential
leads to the decrease of the ground state amplitude
$\varphi$ near the center $x = 0$ of the harmonic potential.

Justification of the asymptotic approximations above for the ground state
of the existence problem (\ref{stationaryGP-last}) appears to be a simple analytical problem
if $W$ is bounded and decaying, while $\eta = \mathcal{O}(\eps^{-1})$
as $\eps \to 0$. We do not include this justification analysis in the present work.

\vspace{0.5cm}

{\bf Acknowledgement.} C.G. is supported by the project ANR-12-MONU-0007 BECASIM.
D.P. is supported by the CNRS Visiting Fellowship. D.P. thanks
members of Institut de Math\'ematiques et de Mod\'elisation, Universit\'e Montpellier
for hospitality and support during his visit (September-November, 2013).

\end{document}